\documentclass[11pt]{article}
\usepackage[T1]{fontenc}
\usepackage[latin1]{inputenc}
\usepackage{amsmath,amssymb}
\usepackage{graphics}
\usepackage{epic}
\usepackage{epsfig}
\usepackage{subfigure}

\newtheorem{thm}{Theorem}[section]
\newtheorem{lem}[thm]{Lemma}
\newtheorem{cor}[thm]{Corollary}
\newtheorem{nota}[thm]{Notation}
\newtheorem{defi}[thm]{Definition}
\newtheorem{prop}[thm]{Proposition}
\newtheorem{rema}[thm]{Remark}

\newenvironment{proof}{\noindent {\bf Proof \phantom{9}}}
{\hfill $\square$ \vspace{0.25cm}}

\numberwithin{equation}{section}
\def\theequation{\thesection.\arabic{equation}}
\def\be{\begin{eqnarray}}
\def\ee{\end{eqnarray}}
\def\numero{\refstepcounter{equation} (\theequation)}

\let\text=\textstyle

\newcommand{\loi}{{\cal L}}

\newcommand{\sm}{{s-}}
\newcommand{\ds}{\displaystyle}
\newcommand{\intot}{\displaystyle \int _0^t }

\newcommand{\indiq}{{\bf 1}}

\newcommand{\e}{{\epsilon}}

\newcommand{\intrd}{\ds{\int_{\rit^d}}}
\newcommand{\rit}{\mathbb{R}}
\newcommand{\nit}{\mathbb{N}}
\newcommand{\dit}{\mathbb{D}}

\title{\bf Individual-based probabilistic models of adaptive evolution and various
  scaling approximations}

\author{Nicolas Champagnat$^{1,2}$, R\'egis Ferri\`ere$^{1,3}$, Sylvie
  M\'el\'eard$^2$}

\date{\today}

\begin{document}

\maketitle

\begin{center}
  \makeatletter\renewcommand{\@makefnmark}{\mbox{$^{\@thefnmark}$}}\makeatother
  \footnotetext[1]{Laboratoire d'Ecologie, Equipe
    Eco-Evolution Math\'ematique, Ecole Normale Sup\'erieure, 46 rue
    d'Ulm, 75230 Paris cedex 05, France }
  \footnotetext[2]{Equipe MODALX
    Universit\'e Paris 10, 200
    avenue de la R\'epublique, 92001 Nanterre Cedex, France}
  \footnotetext[3]{Department of
    Ecology and Evolutionary Biology, University of Arizona, Tucson AZ
    85721, USA }
  \makeatletter\renewcommand{\@makefnmark}{}\makeatother
  \makeatletter\renewcommand{\@makefnmark}{\mbox{$^{\@thefnmark}$}}\makeatother
\end{center}


\begin{abstract}
  We are interested in modelling Darwinian evolution,
  resulting from the interplay of phenotypic variation
  and natural selection through ecological interactions.
  Our models are rooted in the
  microscopic, stochastic description of a population of discrete
  individuals characterized by one or several adaptive traits. The
  population is modelled as a stochastic point process whose generator
  captures the probabilistic dynamics over continuous time of birth,
  mutation, and death, as influenced by each individual's trait
  values, and interactions between individuals. An offspring usually
  inherits the trait values of her progenitor, except when a
  mutation causes the offspring to take an instantaneous
  mutation step at birth to new trait values.  We look for tractable
  large population approximations. By combining various
  scalings on population size, birth and death rates, mutation rate, mutation step, or
  time, a single microscopic model is shown to lead to contrasting
  macroscopic limits,
  of different nature: deterministic, in the form of ordinary, integro-, or partial
  differential equations, or probabilistic, like stochastic partial
  differential equations or superprocesses. In the limit of rare
  mutations, we show that a possible approximation is a jump process,
  justifying rigorously the so-called trait substitution sequence.
  We thus unify different points of
  view concerning mutation-selection evolutionary models.
\end{abstract}

\bigskip



\section{Introduction}
\label{sec:intro}

In this paper, we are interested in modelling  the dynamics of
populations as driven by the interplay of phenotypic variation and
natural selection operating through ecological interactions, i.e.\
Darwinian evolution. The fundamental property of evolving systems
is the propensity of each individual to create and to select the
diversity. This feature requires to focus on the stochastic
dynamics of each individual in the population. The study of such
evolutionary-ecological models is very complicated, and several
approximations have been proposed. Firstly, Bolker and
Pacala~\cite{BP97} and Dieckmann and Law~\cite{DL00} have
introduced the moment equations of the distribution of traits in
the population and studied different moment closure heuristics.
Secondly, various nonlinear macroscopic models
(integro-differential equations, partial differential equations,
superprocesses) have been proposed without microscopic
justification. Finally, the emerging field of adaptive dynamics
have proposed a new class of macroscopic models on the
evolutionary time scale, defined as jump processes and ordinary
differential equations (trait substitution sequences, Metz et
al.~\cite{Mal96}, canonical equation of adaptive dynamics,
Dieckmann and Law~\cite{DL96}). In all these cases and from a
biological point of view, the pathway from microscopic to
macroscopic models deserves a firm mathematical pavement, at least
to clarify the significance of the implicit biological assumptions
underlying the choice of a particular model.

In this work, we unify several macroscopic approximations by
recovering them from a single microscopic model. In particular, we
point out the importance of large population assumptions and that
the nature of the approximation strongly depends on the
combination of various scalings of the biological parameters
(birth and death rates, mutation rate, mutation step and time).

This paper starts (Section~\ref{sec:ppp}) with the microscopic
description of a population of discrete individuals, whose
phenotypes are described by a vector of trait values. The
population is modelled as a stochastic Markov point process whose
generator captures the probabilistic dynamics over continuous time
of birth, mutation and death, as influenced by each individual's
trait values and interactions between individuals. The adaptive
nature of a trait implies that an offspring usually inherits the
trait values of her progenitor, except when a mutation occurs. In
this case, the offspring makes an instantaneous mutation step at
birth to new trait values. We will refer to the state space
parameterized by adaptive traits as the trait space, and will
often (slightly abusively) call trait the actual trait value. This
process is defined as the solution of a stochastic differential
equation driven by point Poisson measures
(Section~\ref{sec:constr}). In Section~\ref{sec:simul}, we give an
algorithmic construction of the population point process and
propose some simulations, for various parameters, of an
asymmetrical example developed in Kisdi~\cite{Ki99}.  Next, we
prove that the point population process is a measure-valued
semimartingale and compute its characteristics
(Section~\ref{secexist}). Then we look for tractable
approximations, following different mathematical paths.  Our first
approach (Section~\ref{sec:moment}) aims at deriving deterministic
equations to describe the moments of trajectories of the point
process, i.e.\ the statistics of a large number of independent
realizations of the process. We explain the difficult hierarchy
between these equations coming from competition kernels and
preventing, even in the simple mean-field case, decorrelations and
tractable moment closure. The alternative approach involves
renormalizations of the point process based on a large population
limit. The measure-valued martingale properties of the
renormalized point process allow us to show that, according to
different scalings of birth, death and mutation rates, one obtains
qualitatively different limiting partial differential equations
and the appearance or not of some demographic stochasticity. We
show in Section~\ref{sec:limit1} that by itself, the
large-population limit leads to a deterministic, nonlinear
integro-dif\-fe\-ren\-tial equation. Then, in
Section~\ref{sec:limit2}, we combine the large-population limit
with an acceleration of birth (hence mutation) and death according
to small mutation steps. That yields either a deterministic
nonlinear reaction-diffusion model, or a stochastic measure-valued
process (depending on the acceleration rate of the birth-and-death
process).  If now this acceleration of birth and death is combined
with a limit of rare mutations, the large-population limit yields
a nonlinear integro-differential equation either deterministic or
stochastic, depending here again on the speed of the scaling of
the birth-and-death process, as described in
Section~\ref{sec:limit3}.

In Section~\ref{sec:AdaDyn}, we model a time scale separation
between ecological events (fast births and deaths) and evolution
(rare mutations), for an initially monomorphic population. The
competition between individuals takes place on the short time
scale. In a large population limit, this leads on the mutation
time scale to a jump process over the trait space, where the
population stays monomorphic at any time.  Thereby we provide a
rigorous justification to the notion of trait substitution
sequence introduced by Metz et al.~\cite{MNG92}.




\section{Population point process}
\label{sec:ppp}

Even if the evolution manifests itself as a global change in the state
of a population, its basic mechanisms, mutation and selection, operate
at the level of individuals. Consequently, we model the evolving
population as a stochastic interacting individual system, where each
individual is characterized by a vector of phenotypic trait values.
The trait space ${\cal X}$ is assumed to be a closed subset of
$\rit^d$, for some $d\geq 1$.

We will denote by $M_F({\cal X})$ the set of finite non-negative
measures on ${\cal X}$. Let also ${\cal M}$ be the subset of
$M_F({\cal X})$ consisting of all finite point measures:
\begin{equation*}
  {\cal M} = \left\{ \sum_{i=1}^n \delta_{x_i} , \; n \geq 0,
    x_1,...,x_n \in {\cal X} \right\}.
\end{equation*}
Here and below, $\delta_x$ denotes the Dirac mass at $x$. For any
$m\in M_F({\cal X})$, any measurable function $f$ on ${\cal X}$, we
set $\left< m, f \right> = \int_{{\cal X}} f dm$.

We aim to study the stochastic process $\nu_t$, taking its values in
${\cal M}$, and describing the distribution of individuals and traits
at time $t$. We define
\begin{equation}
  \label{pop}
  \nu_t = \sum_{i=1}^{I(t)} \delta_{X^i_t},
\end{equation}
$I(t) \in {\mathbb{N}}$ standing for the number of individuals
alive at time $t$, and $X^1_t,...,X^{I(t)}_t$ describing the
individual's traits (in ${\cal X}$).\\

  For a population $\nu=\sum_{i=1}^{I}\delta_{x^i}$, and a trait $x\in
  {\cal X}$, we define
    the birth rate $b(x,V*\nu(x))=b(x,\sum_{i=1}^{I}V(x-x^i))$ and the death rate
      $d(x,U*\nu(x))=d(x,\sum_{i=1}^{I}U(x-x^i))$ of individuals with trait $x$; $V$
    and $U$ denote the interaction kernels affecting respectively
    reproduction and mortality.
 let $\mu(x)$ and $M(x,z)dz$ be respectively the probability that an offspring produced by
    an individual with trait $x$ carries a mutated trait and the law of this mutant
    trait.

Thus, the population evolution can be roughly summarized as
follows. The initial population is characterized by a (possibly
  random) counting measure $\nu_0\in {\cal M}$ at time $0$, and
  any
 individual with trait $x$ at time $t$ has two independent random exponentially
  distributed ``clocks'': a birth clock with parameter
  $b(x,V*\nu_t(x))$, and a death clock with parameter
  $d(x,U*\nu_t(x))$.
If the death clock of an individual rings, this individual dies
  and disappears.
If the birth clock of an individual with trait $x$ rings, this
  individual produces an offspring. With probability $1-\mu(x)$ the
  offspring carries the same trait $x$; with probability $\mu(x)$ the
  trait is mutated.
 If a mutation occurs, the mutated offspring instantly acquires a
  new trait $z$, picked randomly according to the mutation step
  measure $M(x,z)dz$.

Thus we are looking for a ${\cal M}$-valued Markov process
$(\nu_t)_{t\geq 0}$ with infinitesimal generator $L$, defined for real
bounded functions $\phi$ by
\begin{align}
\label{generator}
 L\phi(\nu)  & = \sum_{i=1}^{I}b(x^i,V*\nu(x^i))(1-\mu(x^i))
  (\phi(\nu+\delta_{x^i})-\phi(\nu)) \notag \\
  & +\sum_{i=1}^{I} b(x^i,V*\nu(x^i))\mu(x^i)
 \int_{\cal X}(\phi(\nu+\delta_{z})-\phi(\nu))M(x^i,z)dz\notag \\
  & +\sum_{i=1}^{I}d(x^i,U*\nu(x^i))
  (\phi(\nu-\delta_{x^i})-\phi(\nu)).
\end{align}
The first term of~(\ref{generator}) captures the effect on the
population of birth without mutation; the second term that of
birth with mutation, and the last term that of death. The
density-dependence makes all terms nonlinear.

\subsection{Process construction}
\label{sec:constr}

Let us justify the existence of a Markov process admitting $L$ as
infinitesimal generator. The explicit construction of $(\nu_t)_{t\geq
  0}$ also yields three side benefits: providing a rigorous and
efficient algorithm for numerical simulations (given hereafter),
laying the mathematical basis to derive the moment equations of the
process (Section~\ref{sec:moment}), and establishing a general method
that will be used to derive some large population limits
(Sections~\ref{sec:large-popu} and~\ref{sec:AdaDyn}).

We make the biologically natural assumption that the trait dependency
of birth parameters is ``bounded'', and at most linear for the death
rate. Specifically, we assume

{\bf Assumptions (H):}

There exist constants $\bar{b}$, $\bar{d}$, $\bar{U}$, $\bar{V}$
and $C$ and a probability density function $\bar{M}$ on $\rit^d$
such that for each $\nu=\sum_{i=1}^{I}\delta_{x^i}$ and for
$x,z\in {\cal X}$,
\begin{gather*}
  b(x,V*\nu(x))\leq \bar{b},\quad d(x,U*\nu(x))\leq
  \bar{d}(1+I), \\
  U(x)\leq \bar{U},\quad V(x)\leq \bar{V}, \\
  M(x,z)\leq C\bar{M}(z-x).
\end{gather*}

These assumptions ensure that there exists a constant $\bar{C}$, such
that the total event rate, for a population counting measure
$\nu=\sum_{i=1}^{I}\delta_{x^i}$, obtained as the sum of all event
rates, is bounded by $\ \bar{C}I(1+I)\ $.

Let us now give a pathwise description of the population process
$(\nu_t)_{t\geq 0}$. We introduce the following notation.

\begin{nota}
  \label{defh}
  Let $\mathbb{N}^*=\mathbb{N}\backslash \{0\}$. Let
  $H=(H^1,...,H^k,...): {\cal M} \mapsto
  (\mathbb{R}^d)^{\mathbb{N}^*}$ be defined by
$ H\left(\textstyle\sum_{i=1}^n \delta_{x_i}\right) =
    (x_{\sigma(1)},...,x_{\sigma(n)},0,...,0,...)$,
  where $x_{\sigma(1)}\curlyeqprec ... \curlyeqprec x_{\sigma(n)}$,
  for some arbitrary order $\curlyeqprec$ on $\mathbb{R}^d$ (
  for example the lexicographic order).
\end{nota}
This function $H$ allows us to overcome the following (purely
notational) problem. Choosing a trait uniformly among all traits
in a population $\nu \in {\cal M}$ consists in choosing $i$
uniformly in $\{1,...,\left<
  \nu,1\right>\}$, and then in choosing the individual {\it number}
$i$ (from the arbitrary order point of view). The trait value of
such an individual is thus $H^i(\nu)$.\\

We now introduce the probabilistic objects we will need.

\begin{defi}
  \label{poisson}
  Let $(\Omega, {\cal F}, P)$ be a (sufficiently large) probability
  space. On this space, we consider
  the following four independent random elements:
  \begin{description}
  \item[\textmd{(i)}] a ${\cal M}$-valued random variable $\nu_0$ (the
    initial distribution),
  \item[\textmd{(ii)}] independent Poisson point measures
    $M_1(ds,di,d\theta)$, and $M_3(ds,di,d\theta)$ on $[0,\infty)
    \times \mathbb{N}^*\times\rit^+$, with the same intensity measure
    $\:  ds \left(\sum_{k\geq 1} \delta_k (di) \right) d\theta\:$ (the
    "clonal" birth and the death Poisson measures),
  \item[\textmd{(iii)}] a Poisson point measure
    $M_2(ds,di,dz,d\theta)$ on $[0,\infty) \times \mathbb{N}^* \times
    {\cal X} \times \rit^+$, with intensity measure $\: ds
    \left(\sum_{k\geq 1} \delta_k (di) \right) dz d\theta\:$ (the
    mutation Poisson measure).
  \end{description}
  Let us denote by $( {\cal F}_t)_{t\geq 0}$ the canonical filtration
  generated by these processes.
\end{defi}

We finally define the population process in terms of these stochastic
objects.

\begin{defi}
  \label{dbpe}
  Assume $(H)$. A $( {\cal F}_t)_{t\geq 0}$-adapted stochastic process
  $\nu=(\nu_t)_{t\geq 0}$  is called
  a population process if a.s., for all $t\geq 0$,
  \begin{align}
    \nu_t &= \nu_0 + \int_{[0,t] \times \mathbb{N}^*
      \times\rit^+}\delta_{H^i(\nu_\sm)} \indiq_{\{i \leq \left<
        \nu_\sm,1 \right>\}} \notag \\ &\hskip 2cm \indiq_{\left\{\theta \leq
        b(H^i(\nu_\sm),V*\nu_\sm(H^i(\nu_\sm)))(1-\mu(H^i(\nu_\sm)))\right\}}
    M_1(ds,di,d\theta) \notag \\ &+ \int_{[0,t] \times \mathbb{N}^* \times
      {\cal X}\times\rit^+} \delta_{z } \indiq_{\{i \leq \left<
        \nu_\sm,1 \right>\}} \notag \\ &\hskip 2cm \indiq_{\left\{\theta \leq
        b(H^i(\nu_\sm),V*\nu_\sm(H^i(\nu_\sm)))
        \mu(H^i(\nu_\sm))M(H^i(\nu_\sm),z)
      \right\}} M_2(ds,di,dz,d\theta) \notag \\ &- \int_{[0,t] \times
      \mathbb{N}^* \times\rit^+} \delta_{H^i(\nu_\sm)} \indiq_{\{i \leq
      \left< \nu_\sm,1 \right>\}} \indiq_{\left\{\theta \leq
        d(H^i(\nu_\sm), U*\nu_\sm(H^i(\nu_\sm))) \right\}}
    M_3(ds,di,d\theta) \label{bpe}
  \end{align}
\end{defi}

Let us now show that if $\nu$ solves~(\ref{bpe}), then $\nu$
follows the Markovian dynamics we are interested in.

\begin{prop}
  \label{gi}
  Assume $(H)$ and consider a solution $(\nu_t)_{t\geq 0}$ of
  Eq.~(\ref{bpe}) such that $E(\sup_{t\geq
    T}\langle\nu_t,\mathbf{1}\rangle^2)<+\infty,\ \forall T>0$. Then
  $(\nu_t)_{t\geq 0}$ is a Markov process. Its infinitesimal generator
  $L$ is defined for all bounded and measurable maps $\phi: {\cal
    M}\mapsto \mathbb{R}$, all $\nu \in {\cal M}$,
  by~(\ref{generator}). In particular, the law of $(\nu_t)_{t\geq 0}$
  does not depend on the chosen order $\curlyeqprec$.
\end{prop}

\begin{proof}
  The fact that $(\nu_t)_{t\geq 0}$ is a Markov process is classical.
  Let us now consider a function $\phi$ as in the statement.
   With our notation, $\nu_0 = \sum_{i=1}^{\left< \nu_0, 1\right>}
  \delta_{H^i(\nu_0)}$.  A simple computation, using the
  fact that a.s., $\phi(\nu_t) = \phi(\nu_0) + \sum_{s\leq t}
  (\phi(\nu_\sm + (\nu_s-\nu_\sm))- \phi(\nu_\sm))$, shows that
  \begin{align*}
    \phi(\nu_t) & = \phi(\nu_0)+\int_{[0,t] \times
      \mathbb{N}^* \times\rit^+}\left(\phi(\nu_\sm +
      \delta_{H^i(\nu_\sm)}) - \phi(\nu_\sm) \right)\indiq_{\{i \leq
      \left< \nu_\sm,1 \right>\}} \\ &\hskip 2cm\indiq_{\left\{\theta
        \leq
        b(H^i(\nu_\sm),V*\nu_\sm(H^i(\nu_\sm)))(1-\mu(H^i(\nu_\sm)))\right\}}
    M_1(ds,di,d\theta) \\ &+ \int_{[0,t] \times \mathbb{N}^* \times
      {\cal X}\times\rit^+} \left(\phi(\nu_\sm + \delta_{z}) -
      \phi(\nu_\sm) \right) \indiq_{\{i \leq \left< \nu_\sm,1
      \right>\}} \\ &\hskip 2cm \indiq_{\left\{\theta \leq
        b(H^i(\nu_\sm),V*\nu_\sm(H^i(\nu_\sm)))\mu(H^i(\nu_\sm))
        M(H^i(\nu_\sm),z)
      \right\}} M_2(ds,di,dz,d\theta) \\ &+ \int_{[0,t] \times
      \mathbb{N}^* \times\rit^+} \left(\phi(\nu_\sm -
      \delta_{H^i(\nu_\sm)}) - \phi(\nu_\sm) \right) \indiq_{\{i \leq
      \left< \nu_\sm,1 \right>\}} \\ &\hskip 2cm \indiq_{\left\{\theta
        \leq d(H^i(\nu_\sm), U*\nu_\sm(H^i(\nu_\sm))) \right\}}
    M_3(ds,di,d\theta) .
  \end{align*}
  Taking expectations, we obtain
  \begin{align*}
    E(\phi & (\nu_t))=E(\phi(\nu_0)) \\ &+\intot
    E\Big(\sum_{i=1}^{\left< \nu_s, 1\right>}\bigg\{\left(\phi(\nu_s +
      \delta_{H^i(\nu_s)}) - \phi(\nu_s)
    \right)b(H^i(\nu_s),V*\nu_s(H^i(\nu_s)))(1-\mu(H^i(\nu_s))) \\
    &+\int_{{\cal X}} \left(\phi(\nu_s + \delta_{z}) - \phi(\nu_s)
    \right)b(H^i(\nu_s),V*\nu_s(H^i(\nu_s)))\mu(H^i(\nu_s))M(H^i(\nu_s),z)dz
    \\ &+ \left(\phi(\nu_s - \delta_{H^i(\nu_s)}) - \phi(\nu_s)
    \right)d(H^i(\nu_s), U*\nu_s(H^i(\nu_s)))\bigg\}\Big)ds 
  \end{align*}
  Differentiating this expression at $t=0$ leads to~(\ref{generator}).
\end{proof}

Let us show existence and moment properties for the population
process.

\begin{thm}\label{existence}
  \begin{description}
  \item[\textmd{(i)}] Assume~(H) and that $E \left( \left<\nu_0,1
      \right> \right) <\infty$. Then the process $(\nu_t)_{t\geq 0}$
    defined by Definition \ref{dbpe} is well defined on $\rit_+$.
  \item[\textmd{(ii)}] If furthermore for some $p \geq 1$, $E \left(
      \left<\nu_0,1 \right>^p \right) <\infty$, then for any
    $T<\infty$,
    \begin{equation}
      \label{lp}
      E (\sup_{t\in[0,T]} \left< \nu_t,1\right>^p ) <\infty.
    \end{equation}
  \end{description}
\end{thm}

\begin{proof}
  We first prove~(ii). Consider the process $(\nu_t)_{t\geq 0}$. We
  introduce for each $n$ the stopping time $\tau_n = \inf \left\{ t
    \geq 0 , \; \left<\nu_t,1 \right> \geq n \right\}$. Then a simple
  computation using Assumption~(H) shows that, neglecting the
  non-positive death terms,
  \begin{align*}
    \sup_{s\in[0,t\land \tau_n]} & \left< \nu_s,1\right>^p\leq \left<
      \nu_0,1\right>^p + \int_{[0,t\land \tau_n] \times
      \mathbb{N}^* \times\rit^+} \left((\left< \nu_\sm,1\right> +1)^p-
      \left< \nu_\sm,1\right>^p\right) \indiq_{\{i \leq \left< \nu_\sm,1
      \right>\}} \\ &\hskip 2cm\indiq_{\left\{\theta \leq
        b(H^i(\nu_\sm),V*\nu_\sm(H^i(\nu_\sm)))(1-\mu(H^i(\nu_\sm)))\right\}}
    M_1(ds,di,d\theta) \\ &+ \int_{[0,t] \times \mathbb{N}^* \times
      {\cal X}\times\rit^+} \left((\left< \nu_\sm,1\right> +1)^p- \left<
        \nu_\sm,1\right>^p\right)\indiq_{\{i \leq \left< \nu_\sm,1
      \right>\}} \\ &\hskip 2cm\indiq_{\left\{\theta \leq
        b(H^i(\nu_\sm),V*\nu_\sm(H^i(\nu_\sm)))\mu(H^i(\nu_\sm))
        M(H^i(\nu_\sm),z)\right\}} M_2(ds,di,dz,d\theta).
  \end{align*}
  Using the inequality $(1+x)^p-x^p\leq C_p(1+x^{p-1})$ and taking
  expectations, we thus obtain, the value of $C_p$ changing from line
  to line,
  \begin{align*}
    E(\sup_{s \in [0,t\land \tau_n]}\left< \nu_s,1\right>^p)
     &\leq C_p\left(1 + E \left( \int_0^{t\land \tau_n} \bar{b}
        \left(\left< \nu_\sm,1\right>+ \left< \nu_\sm,1\right>^p
        \right) ds \right)\right) \\ &\leq C_p\left(1+ E\left( \intot
        \left(1+\left< \nu_{s\land \tau_n},1\right>^p \right)
        ds\right)\right).
  \end{align*}
  The Gronwall Lemma allows us to conclude that for any $T<\infty$,
  there exists a constant $C_{p,T}$, not depending on $n$, such that
  \begin{equation}
    \label{xxxx}
    E( \sup_{t \in [0,T\land \tau_n]} \left< \nu_t,1\right>^p)\leq C_{p,T}.
  \end{equation}
  First, we deduce that $\tau_n$ tends a.s. to infinity. Indeed, if
  not, one may find a $T_0<\infty$ such that $\e_{T_0}=P\left( \sup_n
    \tau_n <T_0 \right) >0$. This would imply that $E\left( \sup_{t
      \in [0,T_0\land \tau_n]}\left< \nu_t,1\right>^p \right) \geq
  \e_{T_0} n^p$ for all $n$, which contradicts~(\ref{xxxx}).
We may let $n$ go to infinity in~(\ref{xxxx}) thanks to the Fatou
  Lemma. This leads to~(\ref{lp}).

  Point~(i) is a consequence of point~(ii). Indeed, one builds the
  solution $(\nu_t)_{t\geq 0}$ step by step. One only has to check
  that the sequence of jump instants $T_n$ goes a.s.\ to infinity as
  $n$ tends to infinity. But this follows from~(\ref{lp}) with $p=1$.
\end{proof}

\subsection{Examples and simulations}
\label{sec:simul}

Let us remark that Assumption (H) is satisfied in the case where
\begin{equation}
  b(x,V*\nu(x))=b(x),\quad d(x,U*\nu(x))=d(x)+\alpha(x)\int_{{\cal
  X}}U(x-y)\nu(dy),
\end{equation}
where $b$, $d$ and $\alpha$ are bounded functions.

In the case where moreover, $\mu\equiv 1$, this individual-based
model can also be interpreted as a model of ``spatially structured
population'', where the trait is viewed as a spatial location and
the mutation at each birth event is viewed as dispersal. This kind
of models have been introduced by Bolker and
Pacala~(\cite{BP97,BP99}) and Law et al.~(\cite{LMD03}), and
mathematically studied by Fournier and M\'el\'eard~\cite{FM04}.
The  case $\ U\equiv 1\ $ corresponds to a density-dependence in
the
total population size.\\

We will consider later the particular set of parameters for the
logistic interaction model, taken from Kisdi~\cite{Ki99} and
corresponding to a model of asymmetrical competition:
\begin{gather}
  \bar{\cal X}=[0,4],\quad d(x)=0,\quad\alpha(x)=1,\quad\mu(x)=\mu,
  \notag \\ b(x)=4-x,\quad U(x-y)=\frac{2}{K}\bigg(1-{1\over 1+1,2
  \exp(-4(x-y))}\bigg)\label{eq:ex}
\end{gather}
and $M(x,z)dz$ is a Gaussian law with mean $x$ and variance $\sigma^2$
conditionned to the fact that the mutant stays in $[0,4]$. As we will
see in Section~\ref{sec:large-popu}, the constant $K$ scaling the
strength of competition also scales the population size (when the
initial population size is proportional to $K$). In this model, the
trait $x$ can be interpreted as body size. Equation~(\ref{eq:ex})
means that body size influences the birth rate negatively, and creates
asymmetrical competition reflected in the sigmoid shape of $U$ (being
larger is competitively advantageous).\\

Let us give an algorithmic construction for the population process
(in the general case), simulating the size $I(t)$ of the
population, and the trait vector $\mathbf{X}_t$ of all individuals
alive at time $t$.

At time $t=0$, the initial population $\nu_0$ contains $I(0)$
individuals and the corresponding trait vector is
$\mathbf{X}_0=(X_0^i)_{1\leq i\leq I(0)}$. We introduce the
following sequences of independent random variables, which will
drive the algorithm.
\begin{itemize}
\item The type of birth or death events will be selected according to
  the values of a sequence of random variables $(W_k)_{k\in\nit^*}$
  with uniform law on $[0,1]$.
\item The times at which events may be realized will be described
  using a sequence of random variables $(\tau_k)_{k\in\nit}$ with
  exponential law with parameter $\bar{C}$.
\item The mutation steps will be driven by a sequence of random
variables $(Z_k)_{k\in\nit}$ with law $\bar{M}(z)dz$.
\end{itemize}
We set $T_0=0$ and construct the process inductively for $k\geq 1$ as
follows.

At step $k-1$, the number of individuals is $I_{k-1}$, and the trait
vector of these individuals is $\mathbf{X}_{T_{k-1}}$.

Let $\displaystyle{T_{k}=T_{k-1}+\frac{\tau_k}{I_{k-1}(I_{k-1}+1)}}$.
Notice that $\displaystyle{\frac{\tau_k}{I_{k-1}(I_{k-1}+1)}}$
represents the time between jumps for $I_{k-1}$ individuals, and
$\bar{C}(I_{k-1}+1)$ gives an upper bound on the total event rate for
each individual.

At time $T_k$, one chooses an individual $i_k=i$ uniformly at
random among the $I_{k-1}$ alive in the time interval
$[T_{k-1},T_k)$; its trait is $X^i_{T_{k-1}}$.  (If $I_{k-1}=0$
then $\nu_t=0$ for all $t\geq T_{k-1}$.)
\begin{itemize}
\item If $\displaystyle{0\leq W_k\leq
    \frac{d(X^{i}_{T_{k-1}},\sum_{j=1}^{I_{k-1}}
      U(X^{i}_{T_{k-1}}-X^{j}_{T_{k-1}}))}{\bar{C}(I_{k-1}+1)}=
    W_1^i(\mathbf{X}_{T_{k-1}})}$, then the chosen individual dies,
  and $I_k=I_{k-1}-1$.
\item If $\displaystyle{W_1^i(\mathbf{X}_{T_{k-1}})< W_k \leq
    W_2^i(\mathbf{X}_{T_{k-1}})}$, where
  \begin{equation*}
    W_2^i(\mathbf{X}_{T_{k-1}})=
    W_1^i(\mathbf{X}_{T_{k-1}})+{[1-\mu(X^{i}_{T_{k-1}})]b(X^{i}_{T_{k-1}},
      \sum_{j=1}^{I_{k-1}} V(X^{i}_{T_{k-1}}-X^{j}_{T_{k-1}})) \over
      \bar{C}(I_{k-1}+1)},
  \end{equation*}
  then the chosen individual gives birth to an offspring with
  trait $X^i_{T_{k-1}}$, and $I_k=I_{k-1}+1$.
\item If $\displaystyle{W_2^i(\mathbf{X}_{T_{k-1}})< W_k \leq
    W_3^i(\mathbf{X}_{T_{k-1}}, Z_k)}$, where
  \begin{multline*}
    W_3^i(\mathbf{X}_{T_{k-1}}, Z_k)=W_2^i(\mathbf{X}_{T_{k-1}})+ \\
    {\mu(X^{i}_{T_{k-1}})b(X^{i}_{T_{k-1}}, \sum_{j=1}^{I_{k-1}}
      V(X^{i}_{T_{k-1}}-X^{j}_{T_{k-1}}))M(X^{i}_{T_{k-1}},
      X^{i}_{T_{k-1}}+Z_k)\over \bar{C}\bar{M}(Z_k)(I_{k-1}+1)},
  \end{multline*}
  then the chosen individual gives birth to a mutant offspring with
  trait $X^i_{T_{k-1}}+Z_k$, and $I_k=I_{k-1}+1$.
\item If $W_k>W_3^i(\mathbf{X}_{T_{k-1}}, Z_k)$, nothing happens, and
  $I_k=I_{k-1}$.
\end{itemize}

\bigskip

Then, at any time $t\geq 0$, the number of individuals is defined
by $I(t)=\sum_{k\geq
  0}1_{\{T_k\leq t< T_{k+1}\}}I_k$ and the population process is obtained as $\nu_t=\sum_{k\geq
  0}1_{\{T_k\leq t<
  T_{k+1}\}}\sum_{i=1}^{I_k}\delta_{X^i_{T_k}}$.\\

\begin{figure}
  \centering
  \mbox{\subfigure[$\mu=0.03$, $K=100$, $\sigma=0.1$.]%
{\epsfig{bbllx=0pt,bblly=0pt,bburx=17.99cm,bbury=22.23cm,%
figure=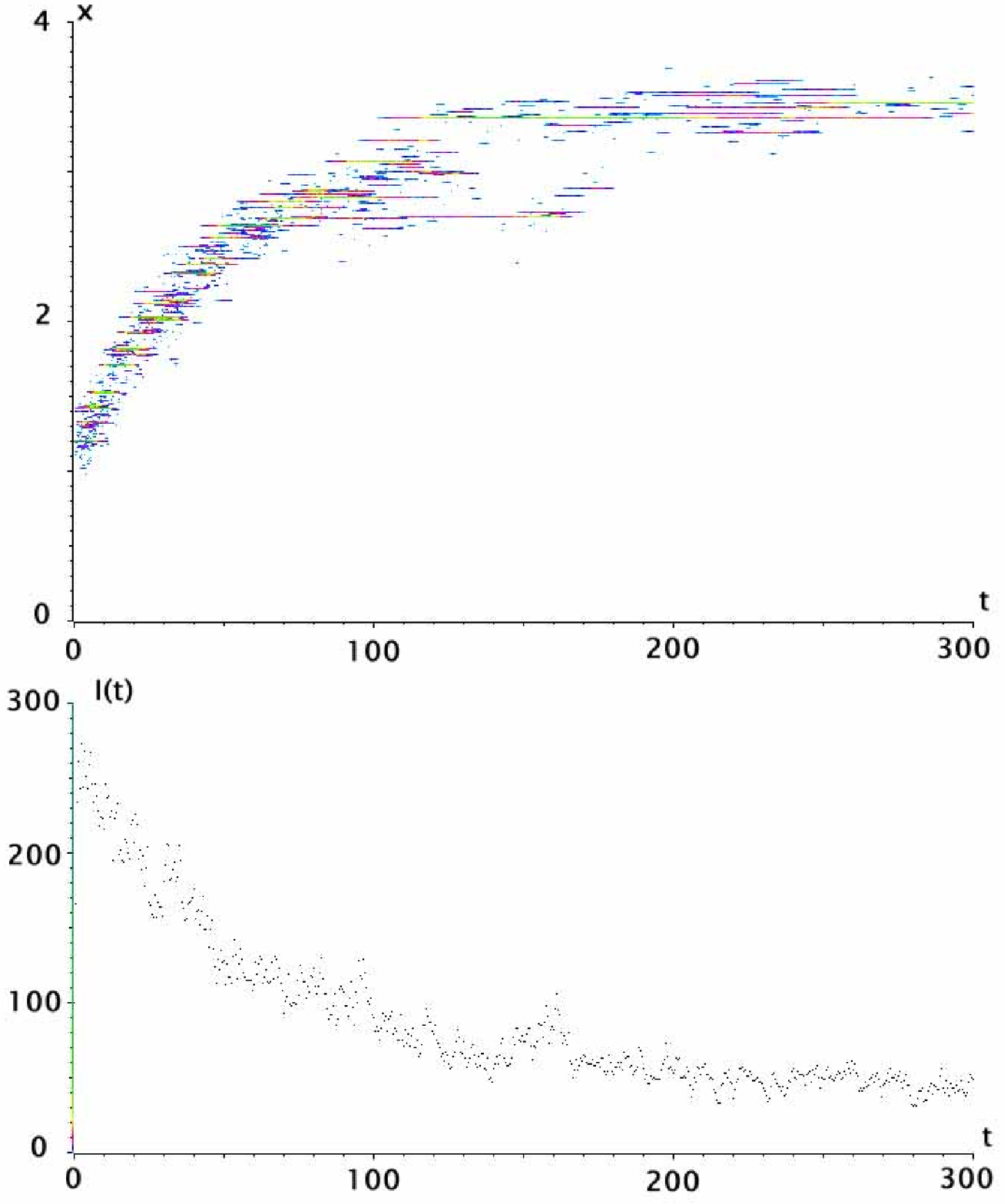, width=.47\textwidth}}\quad
    \subfigure[$\mu=0.03$, $K=3000$, $\sigma=0.1$.]%
{\epsfig{bbllx=0pt,bblly=0pt,bburx=17.99cm,bbury=22.23cm,%
figure=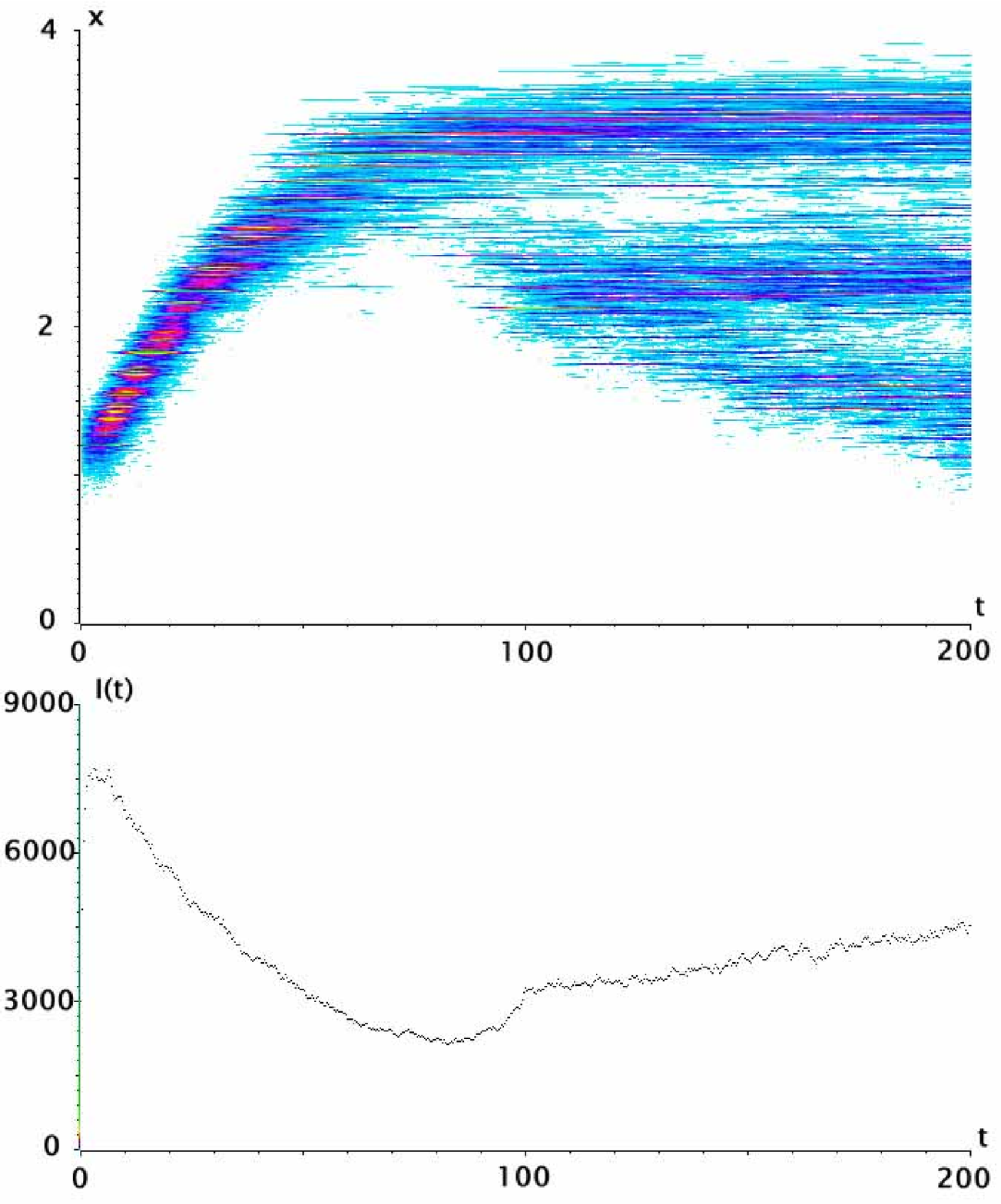, width=.47\textwidth}}} \\
  \mbox{\subfigure[$\mu=0.03$, $K=100000$, $\sigma=0.1$.]%
{\epsfig{bbllx=0pt,bblly=0pt,bburx=17.99cm,bbury=22.23cm,%
figure=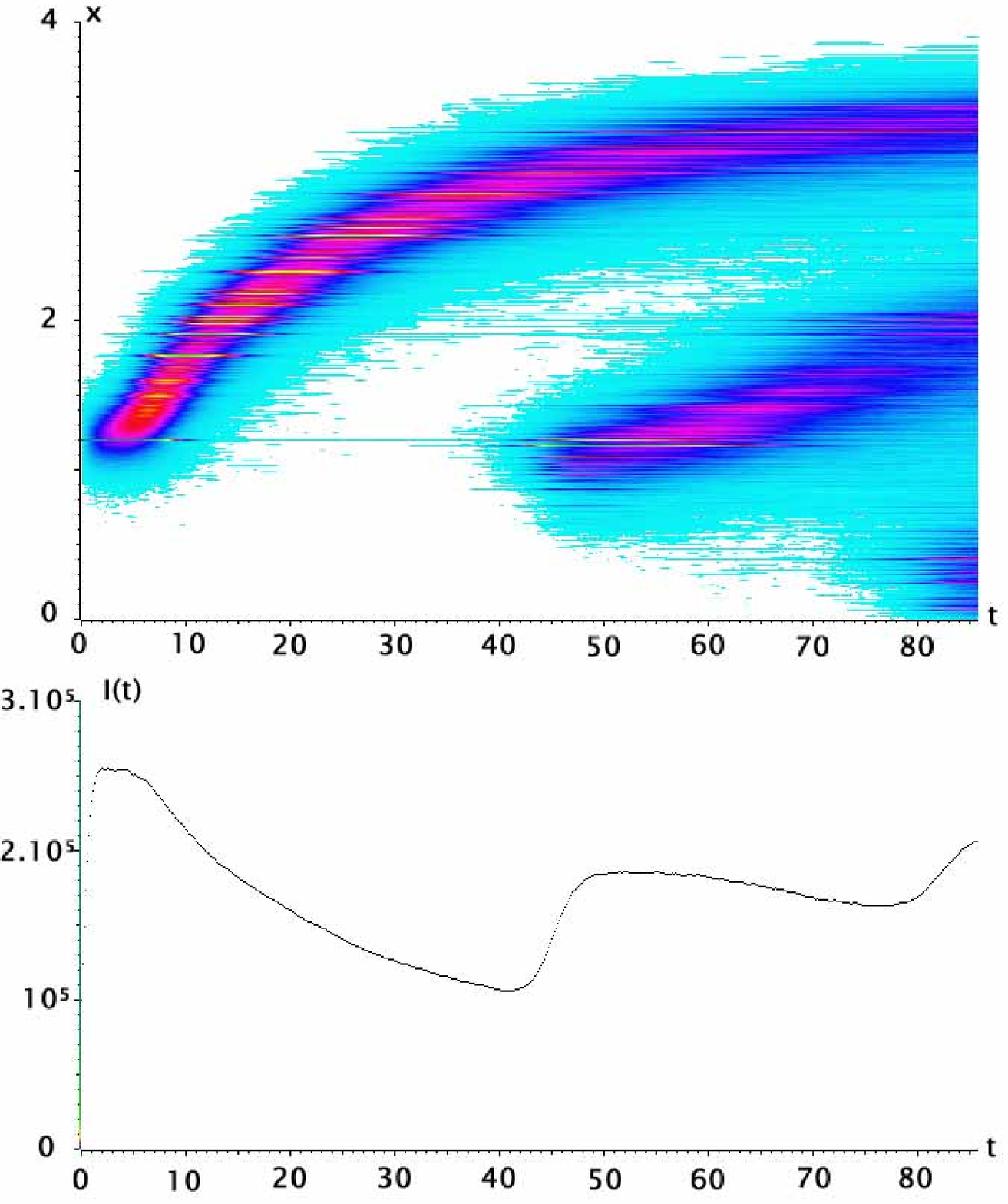, width=.47\textwidth}}\quad
    \subfigure[$\mu=0.00001$, $K=3000$, $\sigma=0.1$.]%
{\epsfig{bbllx=0pt,bblly=0pt,bburx=17.99cm,bbury=22.23cm,%
figure=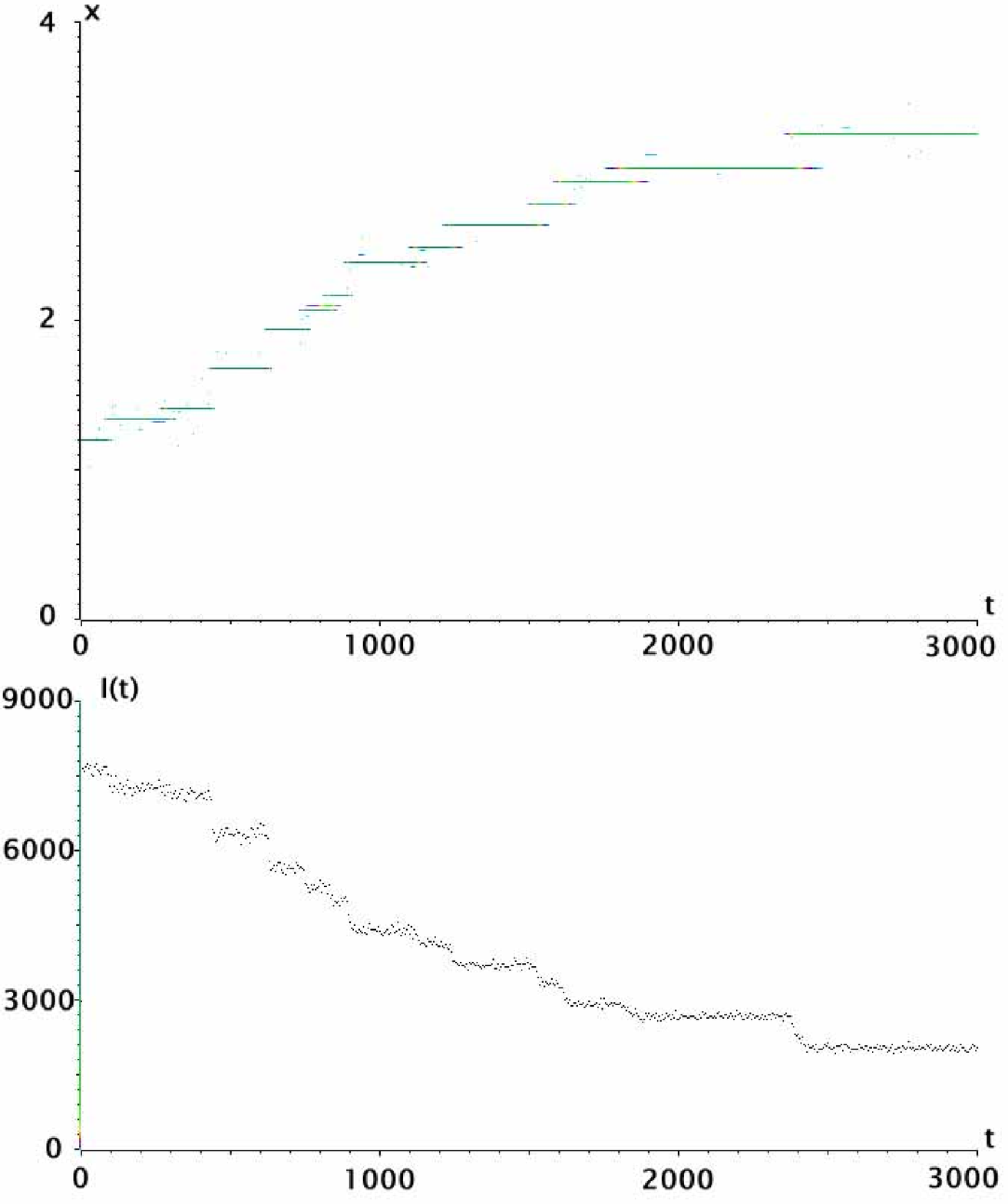, width=.47\textwidth}}}
  \caption{Numerical simulations of trait distributions (upper panels,
    darker is higher frequency) and population size (lower panels).
    The initial population is monomorphic with trait value $1.2$ and
    contains $K$ individuals. (a--c) Qualitative effect of increasing
    system size (measured by parameter $K$). (d) Large system size and
    very small mutation probability ($\mu$).}
  \label{fig:IBM-1}
\end{figure}

The simulation of Kisdi's example~(\ref{eq:ex}) can be carried out
following this algorithm. We can show a very wide variety of
qualitative behavior according to the value of the parameters
$\sigma$, $\mu$ and $K$.

In the following figures, the upper part gives the distribution of
the traits in the population at any time, using a grey scale code
for the number of individuals holding a given trait. The lower
part of the simulation represents the dynamics of the total size
$I(t)$ of the population.

These simulations will serve to illustrate the different mathematical
scalings described in Sections~\ref{sec:large-popu}
and~\ref{sec:AdaDyn}. Let us observe for the moment the qualitative
differences between the cases where $K$ is large
(Fig.~\ref{fig:IBM-1}~(c)), in which a wide population density evolves
regularly (see Section~\ref{sec:limit1}) and where $\mu$ is small
(Fig.~\ref{fig:IBM-1}~(d)), in which the population trait evolves
according to a jump process (see Section~\ref{sec:TSS}).

\begin{figure}
  \centering
  \mbox{\subfigure[$\mu=0.3$, $K=10000$, $\sigma=0.3/K^{\eta/2}$, $\eta=0.5$.]%
{\epsfig{bbllx=0pt,bblly=0pt,bburx=17.99cm,bbury=22.23cm,%
figure=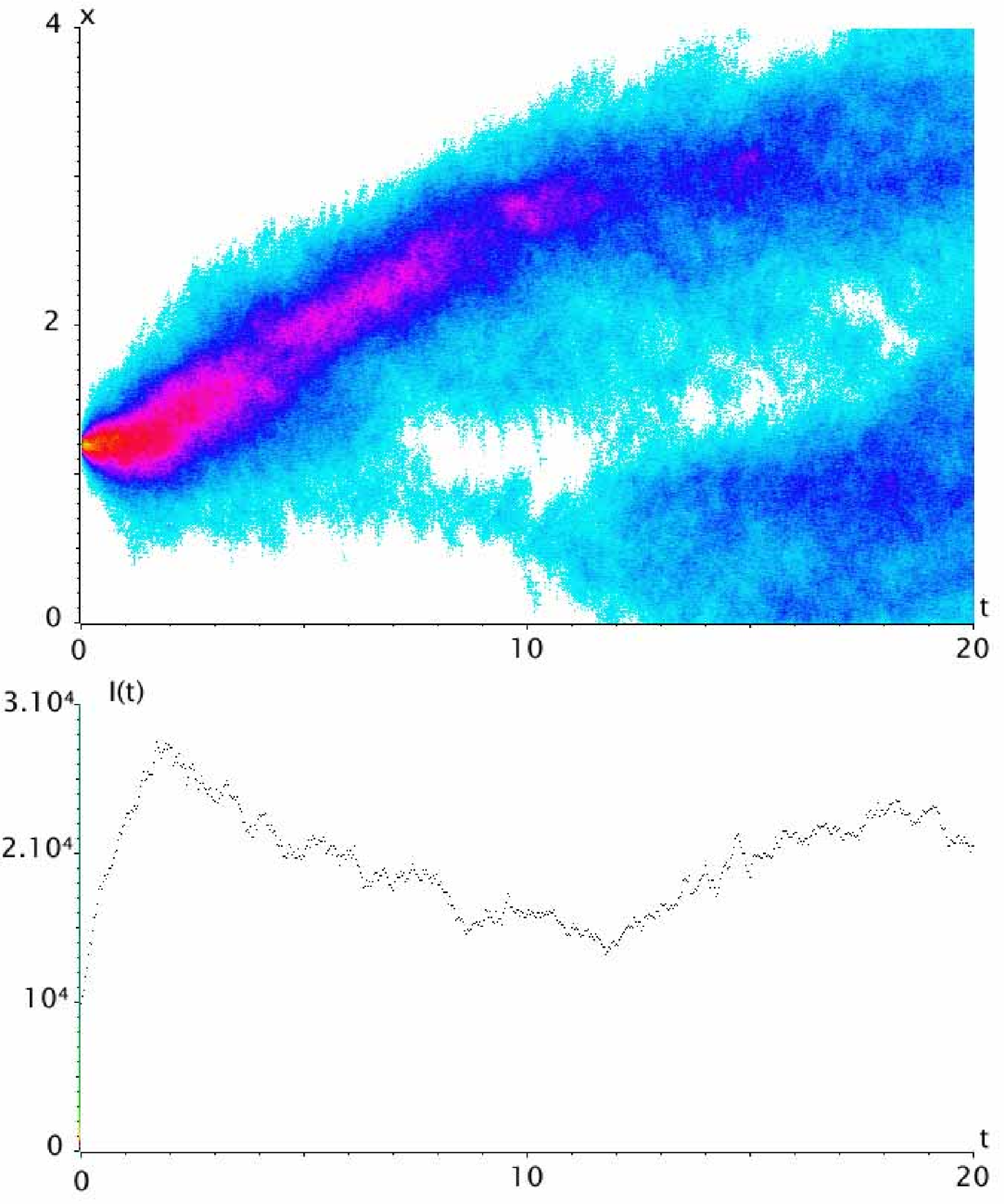, width=.47\textwidth}}\quad
    \subfigure[$\mu=0.1/K^{\eta}$, $K=10000$, $\sigma=0.1$, $\eta=0.5$.]%
{\epsfig{bbllx=0pt,bblly=0pt,bburx=17.99cm,bbury=22.23cm,%
figure=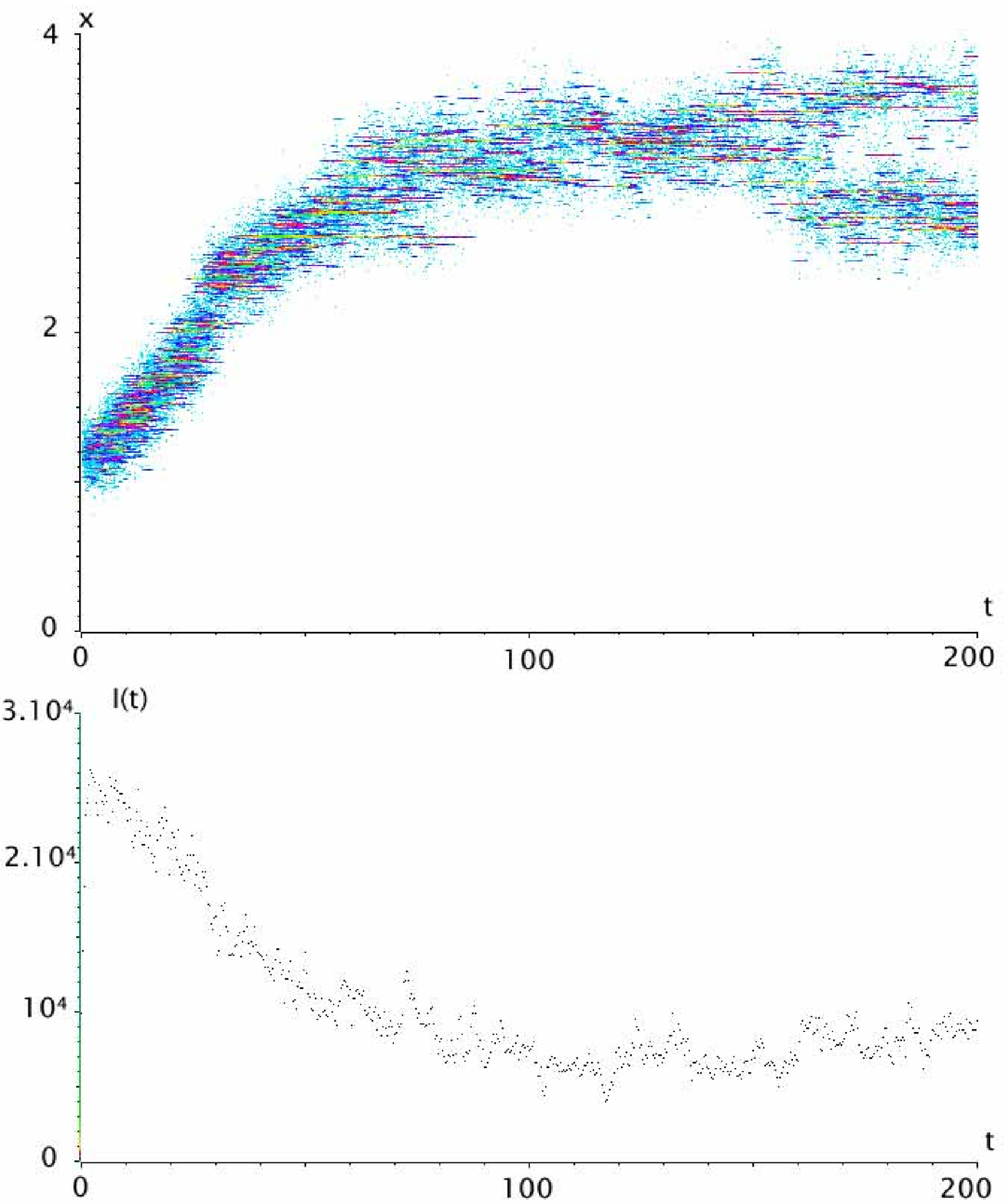, width=.47\textwidth}}} \\
  \mbox{\subfigure[$\mu=0.3$, $K=10000$, $\sigma=0.3/K^{\eta/2}$, $\eta=1$.]%
{\epsfig{bbllx=0pt,bblly=0pt,bburx=17.99cm,bbury=22.23cm,%
figure=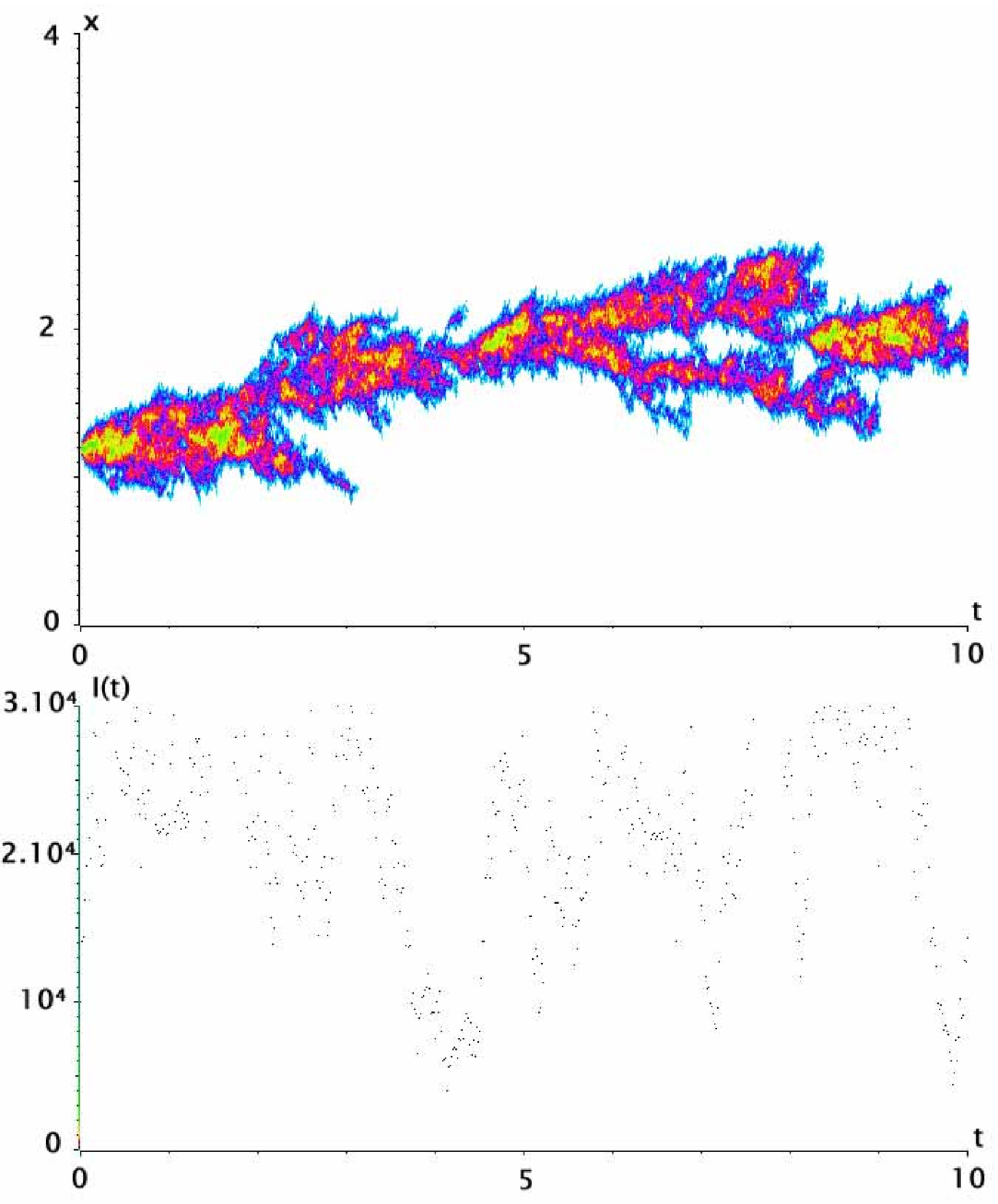, width=.47\textwidth}}\quad
    \subfigure[$\mu=0.3$, $K=10000$, $\sigma=0.3/K^{\eta/2}$, $\eta=1$.]%
{\epsfig{bbllx=0pt,bblly=0pt,bburx=17.99cm,bbury=22.23cm,%
figure=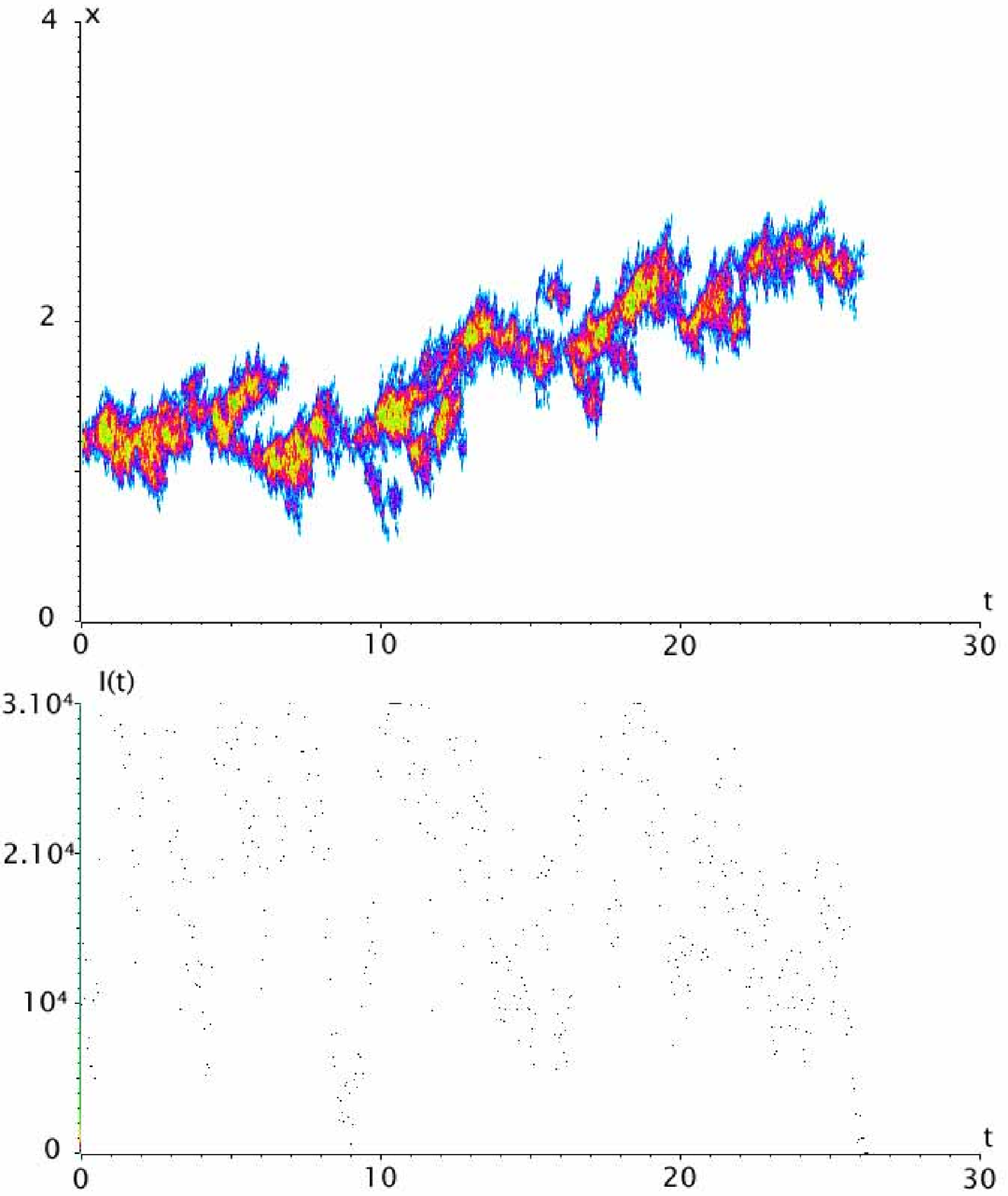, width=.47\textwidth}}}
  \caption{Numerical simulations of trait distribution (upper panels,
    darker is higher frequency) and population size (lower panels) for
    accelerated birth and death and concurrently increased system
    size. Parameter $\eta$ (between 0 and 1) relates the acceleration
    of demographic turnover and the increase of system size. (a)
    Rescaling mutation step. (b) Rescaling mutation probability.
    (c--d) Rescaling mutation step in the limit case $\eta=1$; two
    samples for the same population. The initial population is
    monomorphic with trait value $1.2$ and contains $K$ individuals.}
  \label{fig:IBM-2}
\end{figure}

The simulations of Fig.~\ref{fig:IBM-2} involve an acceleration of the
birth and death processes (see Section~\ref{sec:accel}) as
\begin{equation*}
  b(x,\zeta)=K^{\eta}+b(x)\quad\mbox{and}\quad
  d(x,\zeta)=K^{\eta}+d(x)+\alpha(x)\zeta.
\end{equation*}
There is a noticeable qualitative difference between
Fig.~\ref{fig:IBM-2}~(a) and~(b), where $\eta=1/2$, and
Fig.~\ref{fig:IBM-2}~(c) and~(d), where $\eta=1$. In the latter, we
observe strong fluctuations in the population size and a finely
branched structure of the evolutionnary pattern, revealing a new form
of stochasticity in the large population approximation.


\subsection{Martingale Properties}\label{secexist}

We finally give some martingale properties of the process
$(\nu_t)_{t\geq 0}$, which  are the key point of our approach.

\begin{thm}
  \label{martingales}
  Assume $(H)$, and that for some $p\geq 2$, $E \left( \left<
      \nu_0,1 \right>^p \right) <\infty$.
  \begin{description}
  \item[\textmd{(i)}] For all measurable functions $\phi$ from ${\cal
      M}$ into $\mathbb{R}$ such that for some constant $C$, for all
    $\nu \in {\cal M}$, $\vert \phi (\nu) \vert + \vert L \phi (\nu)
    \vert \leq C (1+\left< \nu,1 \right>^p)$, the process
    \begin{equation}
      \label{pbm2}
      \phi(\nu_t) - \phi(\nu_0) - \intot  L\phi (\nu_s) ds
    \end{equation}
    is a c\`adl\`ag $({\cal F}_t)_{t\geq 0}$-martingale starting from
    $0$.
  \item[\textmd{(ii)}] Point~(i) applies to any function $\phi(\nu)=
    \left< \nu, f \right>^q$, with $0\leq q\leq p-1$ and with $f$
    bounded and measurable on ${\cal X}$.
  \item[\textmd{(iii)}] For such a function $f$, the process
    \begin{align}
      M^f_t &= \langle \nu_t, f\rangle - \langle\nu_0,f\rangle -
      \intot \int_{\cal X}\bigg\{\bigg((1-\mu(x))b(x,V*\nu_s(x))
      -d(x,U*\nu_s(x))\bigg)f(x) \notag \\
      &+\mu(x)b(x,V*\nu_s(x))\int_{\cal X}f(z)M(x,z)dz \bigg\}
      \nu_s(dx) ds \label{eq:mart}
    \end{align}
    is a c\`adl\`ag square integrable martingale starting from $0$
    with quadratic variation
    \begin{align}
      \langle M^f\rangle_t &= \intot \int_{\cal
        X}\bigg\{\bigg((1-\mu(x))b(x,V*\nu_s(x))-d(x,U*\nu_s(x))\bigg)f^2(x)
      \notag \\ &+\mu(x)b(x,V*\nu_s(x))\int_{\cal X}f^2(z)M(x,z)dz
      \bigg\} \nu_s(dx) ds. \label{qv1}
    \end{align}
  \end{description}
\end{thm}

\begin{proof}
  First of all, note that point~(i) is immediate thanks to
  Proposition~\ref{gi} and~(\ref{lp}). Point~(ii) follows from a
  straightforward computation using~(\ref{generator}). To prove~(iii),
  we first assume that $E \left( \left< \nu_0,1 \right>^3 \right)
  <\infty$. We apply~(i) with $\phi(\nu) = \langle \nu, f\rangle$.
  This yields that $M^f$ is a martingale. To compute its bracket, we
  first apply~(i) with $\phi(\nu) = \left< \nu, f\right>^2$ and obtain
  that
  \begin{align}
    \langle \nu_t, f\rangle^2 &- \langle\nu_0,f\rangle^2 - \intot
    \int_{\cal X}\bigg\{\bigg((1-\mu(x))b(x,V*\nu_s(x))(f^2(x) + 2
    f(x)\left< \nu_s, f \right>) \notag \\
    &\hskip 3cm+d(x,U*\nu_s(x))(f^2(x) - 2 f(x)\left< \nu_s, f
    \right>)\bigg) \notag \\ &+\mu(x)b(x,V*\nu_s(x))\int_{\cal X}(f^2(z)
    + 2 f(z) \left< \nu_s, f \right>)M(x,z)dz \bigg\} \nu_s(dx) ds \label{cc1}
  \end{align}
  is a martingale. In another hand, we apply the It\^o formula to
  compute $\left< \nu_t, f\right>^2$ from~(\ref{eq:mart}). We deduce
  that
  \begin{align}
    \langle \nu_t, f\rangle^2 &- \langle\nu_0,f\rangle^2 -\intot
    2\left< \nu_s, f \right>\int_{\cal
      X}\bigg\{\bigg((1-\mu(x))b(x,V*\nu_s(x))-d(x,U*\nu_s(x))\bigg)f(x)
    \notag \\ &+\mu(x)b(x,V*\nu_s(x))\int_{\cal X} f(z) M(x,z)dz
    \bigg\} \nu_s(dx) ds - \langle M^f \rangle_t \label{cc2}
  \end{align}
  is a martingale. Comparing~(\ref{cc1}) and~(\ref{cc2}) leads
  to~(\ref{qv1}). The extension to the case where only $E \left(
    \left< \nu_0,1 \right>^2\right)<\infty$ is straightforward, since
  even in this case, $E(\langle M^f \rangle_t)<\infty$ thanks
  to~(\ref{lp}) with $p=2$.
\end{proof}

\section{Moment equations}
\label{sec:moment}

Moment equations have been proposed by Bolker and
Pacala~(\cite{BP97,BP99}) and Dieckmann and Law~(\cite{DL00}) as
handy analytical models for spatially structured populations.

The philosophy of moment equations is germane to the principle of
Monte-Carlo methods: computing the mean path of the point process
from a large number of independent realizations. (Another
approach, as we shall see in Section~\ref{sec:large-popu}, is to
model the behavior of a single trajectory when it is the initial
number of individuals which is made large).

Let us define the deterministic measure $E(\nu)$ associated with a
random measure $\nu$ by $\displaystyle{\int_{\cal
    X}\varphi(x)E(\nu)(dx)=E(\int_{\cal X}\varphi(x)\nu(dx))}$. Taking
expectations in~(\ref{eq:mart}), we obtain some formula for
$\int_{\cal X}\varphi(x)E(\nu)(dx)$ involving the expectations of
integrals with respect to $\nu(dx)$ or to $\nu(dx)\nu(dy)$.
Nevertheless, this equation is very intricate and presents an
unresolved hierarchy of nonlinearities.  Writing an equation for
$E(\nu(dx)\nu(dy))$ could be possible but will involve integrals
with respect to $\nu(dx)\nu(dy)\nu(dz)$ and so on. Whether this
approach may eventually help describe the population dynamics in
the trait space is still unclear.

Let us consider the case of spatially structured population (see
Section~\ref{sec:simul}) where $d(x,\zeta)=d(x)+\alpha(x)\zeta$, $\:
b(x,\zeta)=b(x)$ and $\mu(x)=1$. Let $N(t)=E(I(t))$ where $I(t)$ is
the number of individuals at time $t$. Taking expectations on
(\ref{eq:mart}) with $\varphi\equiv 1$ yields:
\begin{align}
  \label{eq:evpop}
  N(t)\! =\! N(0)\! +\! \int_0^t\! E\left(\int_{\cal X}\!
    (b(x)-d(x))\nu_s(dx)-\!\int_{{\cal X}\times{\cal X}}\!
    \alpha(x)U(x-y)\nu_s(dx)\nu_s(dy)\right)\! ds.
\end{align}
In the specific case where $b$, $d$ and $\alpha$ are independent
of (the spatial location) $x$, (cf.~\cite{LMD03}),
(\ref{eq:evpop}) recasts into
\begin{equation}
  \label{eq:evpop1}
  \dot{N}=(b-d)N-\alpha E\left(\int_{{\cal X}\times{\cal X}}
   U(x-y)\nu_t(dx)\nu_t(dy)\right).
\end{equation}
Even in the specific mean-field case where $\: U=1\:$, we get
\begin{align}
  \label{eq:meanfield}
  \dot{N}=(b-d)N-\alpha E\left(\int_{{\cal X}\times{\cal X}}
   \nu_t(dx)\nu_t(dy)\right).
\end{align}
The quadratic term corresponding to spatial correlations can not
be simplified and (\ref{eq:meanfield}) allows us to precisely
identify the mathematical issues raised by the problem of moment
closure. In Section~\ref{sec:limit1}, we will see that one needs
the additional large population hypothesis to decorrelate the
quadratic term and to refind the well-known logistic equation.
\medskip

Nevertheless, even if we are not able to produce a closed equation
satisfied by $E(\nu)$, we are able to show, in the general case,
the following qualitative important property concerning the
absolute continuity of the expectation of $\nu_t$.

\begin{prop}
  \label{densBP}
  Assume~(H), that $E(\left< \nu_0,1\right> )<\infty$ and that
  $E(\nu_0)$ is absolutely continuous with respect to the Lebesgue
  measure.  Then for all $t\geq 0$, $E(\nu_t)$ is absolutely continuous
  with respect to the Lebesgue
  measure.
\end{prop}

\begin{rema}
  This implies in particular that, when the initial trait distribution
  $E(\nu_0)$ has no singularity w.r.t.\ the Lebesgue measure, these
  singularities, such as Dirac masses, can only appear in the limit of
  infinite time.
\end{rema}

\begin{proof}
  Consider a Borel set $A$ of $\mathbb{R}^d$ with Lebesgue measure
  zero. Consider also, for each $n\geq 1$, the stopping time $\tau_n =
  \inf \left\{ t \geq 0 , \; \left<\nu_t,1 \right> \geq n \right\}$. A
  simple computation allows us to obtain, for all $t\geq 0$, all $n
  \geq 1$,
  \begin{align*}
    E \left( \left<\nu_{t\land \tau_n}, \indiq_A \right> \right) &\leq
    E(\langle \nu_0,\indiq_A\rangle)+ \bar{b} \ E \bigg(\int_0^{t\land
      \tau_n} \int_{\cal X} \indiq_A(x)\nu_s(dx) ds \bigg) \\
    &+\bar{b} \ E\bigg(\int_0^{t\land \tau_n} \int_{\cal X}
    \left(\int_{\cal X} \indiq_A(z) M(x,z) dz \right)\nu_s(dx)ds
    \bigg).
  \end{align*}
  By assumption, the first term on the RHS is zero. The third term is
  also zero, since for any $x \in {\cal X}$, $\int_{\cal X}
  \indiq_A(z)M(x,z)dz=0$. By Gronwall's lemma, we conclude that for
  each $n$, $E(\langle \nu_{t\land \tau_n},\indiq_A\rangle)$ is zero.
  Thanks to~(\ref{lp}) with $p=1$, $\tau_n$ a.s.\ grows to infinity
  with $n$, which concludes the proof.
\end{proof}

\section{Large-population renormalizations of the individual-based
  process}
\label{sec:large-popu}

The moment equation approach outlined above is based on the idea of
averaging a large number of independent realizations of the population
process initiated with a finite number of individuals.  If $K$ scales
the initial number of individuals, the alternative approach consists
in studying the exact process by letting that system size become very
large and making some appropriate renormalizations.  Several types of
approximations can then be derived, depending on the renormalization
of the process.

For any $K$, let the set of parameters $U_K$, $V_K$, $b_K$, $d_K$,
$M_K$, $\mu_K$ satisfy the Assumption~(H). Let $\nu^K_t$ be the
counting measure of the population at time $t$. We define the
measure-valued Markov process $(X^K_t)_{t\geq 0}$ by
\begin{equation*}
  X^K_t=\frac{1}{K}\nu^K_t.
\end{equation*}

As the system size $K$ goes to infinity, we need to assume the

{\bf Assumption (H1):} The parameters $U_K$, $V_K$, $b_K$, $d_K$,
$M_K$ and $\mu_K$ are all continuous, $\zeta\mapsto b(x,\zeta)$ and
$\zeta\mapsto d(x,\zeta)$ are Lipschitz for any $x\in{\cal X}$, and
\begin{equation*}
  U_K(x)=U(x)/K,\quad V_K(x)=V(x)/K.
\end{equation*}

A biological interpretation of this renormalization is that larger
systems are made up of smaller individuals, which may be a
consequence of a fixed amount of available resources to be
partitioned among individuals.  Thus, the biomass of each
interacting individual scales as $1/K$, which may imply that the
interaction effect of the global population on a focal individual
is of order $1$. Parameter $K$ may also be interpreted as scaling
the resources available, so that the renormalization of $U_K$ and
$V_K$ reflects the decrease of competition for resources.

The generator $\tilde{L}^K$ of $(\nu^K_t)_{t\geq 0}$ is given
by~(\ref{generator}), with parameters $U_K$, $V_K$, $b_K$, $d_K$,
$M_K$, $\mu_K$. The generator $L^K$ of $(X^K_t)_{t\geq 0}$ is obtained
by writing, for any measurable function $\phi$ from $M_F({\cal X})$
into $\rit$ and any $\nu\in M_F({\cal X})$,
\begin{equation*}
  L^K\phi(\nu)=\partial_tE_{\nu}(\phi(X^K_t))_{t=0}= \partial_t
  E_{K\nu}(\phi(\nu^K_t/K))_{t=0} = \tilde{L}^K\phi^K(K\nu)
\end{equation*}
where $\phi^K(\mu)=\phi(\mu/K)$. Then we get
\begin{align}
  \label{eq:def-LK}
  L^K\phi(\nu)&=K\int_{\cal X} b_K(x,V*\nu(x))(1-\mu_K(x))(\phi(\nu+{1\over
    K}\delta_x)-\phi(\nu)) \nu(dx) \notag \\ &+K\int_{\cal X}\int_{\cal X}
    b_K(x,V*\nu(x))\mu_K(x)
    (\phi(\nu+{1\over
    K}\delta_{z})-\phi(\nu))M_K(x,z)dz \nu(dx)\notag \\ &+K\int_{\cal X}
  d_K(x,U*\nu(x))(\phi(\nu-{1\over
    K}\delta_x)-\phi(\nu)) \nu(dx).
\end{align}

By a similar proof as the one of Section~\ref{secexist}, we may
summarize the moment and martingale properties of $X^K$.
\begin{prop}
  \label{XK}
  Assume that for some $p\geq 2$,  $E(\langle X^K_0,1\rangle^p)<+\infty$.
  \begin{description}
  \item[\textmd{(1)}] For any $T>0$, $E(\sup_{t\in[0,T]}\langle X^K_t,1\rangle^p)<+\infty$.
  \item[\textmd{(2)}] For any bounded and measurable functions $\phi$
    on $M_F$ such that $|\phi(\nu)|+|L^K\phi(\nu)|\leq
    C(1+<\nu,1>^p)$, the process
    $\ \phi(X^K_t)-\phi(X^K_0)-\int_0^tL^K\phi(X^K_s)ds\ $
    is a c\`adl\`ag martingale.
  \item[\textmd{(3)}] For each measurable bounded function $f$, the
    process
    \begin{align*}
      m^{K,f}_t & = \langle X^K_t,f\rangle
      -\langle X^K_0, f\rangle \\ &
      -\int_0^t\int_{\cal X}(b_K(x,V* X^K_s(x))-d_K(x,U*X^K_s(x)))f(x)
      X^K_s(dx)ds \\ &-\int_0^t\int_{\cal X}
      \mu_K(x)b_K(x,V*X^K_s(x)
      \bigg(\int_{\cal X}f(z)M_K(x,z)dz-f(x)\bigg)X^K_s(dx)ds
    \end{align*}
    is a square integrable martingale with quadratic variation
  \end{description}
  \vspace{-0.8cm}
  \begin{multline}
    \label{varquad}
    \langle m^{K,f}\rangle_t={1\over
      K}\bigg\{\int_0^t\int_{\cal X}\mu_K(x)b_K(x,V*X^K_s(x))
    \bigg(\int_{\cal X}f^2(z)M_K(x,z)dz- f^2(x)\bigg)X^K_s(dx)ds \\
    +\int_0^t\int_{\cal X}( b_K(x,V*X^K_s(x))+d_K(x,U*
    X^K_s(x)))f^2(x)X^K_s(dx)ds\bigg\}
  \end{multline}
\end{prop}
The search of tractable limits for the semimartingales $\langle
X^K,f\rangle$ yields the different choices of scalings of the
parameters developed in this section.  In particular, we obtain
the deterministic or stochastic nature of the approximation by
studying the quadratic variation of the martingale term, given
in~(\ref{varquad}).

\subsection{Large-population limit}
\label{sec:limit1}

We assume here that $b_K=b$, $d_K=d$, $\mu_K=\mu$, $M_K=M$.

\begin{thm}
  \label{largepop}
  Assume Assumptions~(H) and~(H1). Assume moreover that the initial
  conditions $X^K_0$ converge in law and for the weak topology on
  $M_F({\cal X})$ as $K$ increases, to a finite deterministic measure
  $\xi_0$, and that $\ \sup_KE(\langle X^K_0,1\rangle^3)<+\infty$.

  Then for any $T>0$, the process $(X^K_t)_{t\geq 0}$ converges in
  law, in the Skorohod space $\dit([0,T],M_F({\cal X}))$, as $K$ goes
  to infinity, to the unique deterministic continuous function $\xi\in
  C([0,T],M_F({\cal X}))$ satisfying  for any bounded $f:{\cal
  X}\rightarrow\rit$
  \begin{align}
    \label{eq:limit1}
    \langle\xi_t,f\rangle & =\langle\xi_0,f\rangle
    +\int_0^t\int_{\cal
      X}f(x)[(1-\mu(x))b(x,V*\xi_s(x))-d((x,U*\xi_s(x))]\xi_s(dx)ds
    \notag \\ & +\int_0^t\int_{\cal X}\mu(x)b(x,V*\xi_s(x))
    \left(\int_{\cal X}f(z)M(x,z)dz\right)\xi_s(dx)ds
  \end{align}
\end{thm}
The proof of Theorem~\ref{largepop} is let to the reader. It can be
adapted from the proofs of Theorem~\ref{readif} and~\ref{readifstoch}
below, or obtained as a generalization of Theorem~5.3 in~\cite{FM04}.
This result is illustrated by the simulations of
Figs.~\ref{fig:IBM-1}~(a)--(c).\bigskip

\noindent {\bf Main Examples:}
\begin{description}
\item[(1) A density case.]  Following similar arguments as in the
  proof of Proposition~\ref{densBP}, one shows that if the initial
  condition $\xi_0$ has a density w.r.t.\  Lebesgue measure, then
  the same property holds for the finite measure $\xi_t$,
   which is then solution of the functional
  equation:
  \begin{align}
    \label{eq:EID}
    \partial_t\xi_t(x) & =\left[(1-\mu(x))b(x,V*\xi_t(x))
      -d(x,U*\xi_t(x))\right]\xi_t(x) \notag \\
    & +\int_{\rit^d}M(y,x) \mu(y) b(y,V*\xi_t(y))\xi_t(y)dy
  \end{align}
  for all $x\in{\cal X}$ and $t\geq 0$. Desvillettes et
  al.~\cite{DPF04} suggest to refer to $\xi_t$ as the population
  number density; then the quantity $n_t=\int_{\cal X}\xi_t(x)dx$
  can be interpreted as the total population density over the
  whole trait space.
\item[(2) The mean field case.]  As for moment equations (cf.\
  Section~\ref{sec:moment}), the case of spatially structured
  populations with constant rates $b$, $d$, $\alpha$ is meaningful.
  In this context,~(\ref{eq:EID}) leads to the following equation on
  $\: n_t$:
  \begin{equation}
    \label{eq:phi=1}
    \partial_t n_t= (b-d)n_t -
    \alpha\int_{{\cal X}\times {\cal X}} U(x-y)\xi_t(dx)\xi_t(dy).
  \end{equation}
  With the assumption $U\equiv 1$, we recover the classical mean-field
  logistic equation of population growth:
  \begin{equation*}
    \partial_t n_t= (b-d)n_t - \alpha n_t^2.
  \end{equation*}
  Comparing~(\ref{eq:phi=1}) with the first-moment
  equation~(\ref{eq:meanfield}) obtained previously stresses out the
  ``decorrelative'' effect of the large system size renormalization (only in case $U\equiv
  1$).
  In~(\ref{eq:meanfield}), the correction term capturing the effect of
  spatial correlations in the population remains, even if one
  assumes $U\equiv 1$.
\item[(3) Monomorphic and dimorphic cases without mutation.]  We
  assume here that the population evolves without mutation (parameter
  $\mu=0$); then the population traits are the initial ones.

  {\bf (a) Monomorphic case:} only trait $x$ is present in the population at
  time $t=0$.  Thus, we can write $X^K_0=n^K_0(x)\delta_x$, and then
  $X^K_t=n^K_t(x)\delta_x$ for any time $t$. Theorem~\ref{largepop}
  recasts in this case into $n^K_t(x)\rightarrow n_t(x)$ with
  $\xi_t=n_t(x)\delta_x$, and~(\ref{eq:limit1}) writes
  \begin{equation}
    \label{eq:monomorph}
    \frac{d}{dt}n_t(x)=n_t(x)\big(b(x,V(0)n_t(x))-d(x,U(0)n_t(x))\big),
  \end{equation}

  {\bf (b) Dimorphic case:} when the population contains two traits $x$ and
  $y$, i.e.\ when $X^K_0=n^K_0(x)\delta_x+n^K_0(y)\delta_y$, we can
  define in a similar way $n_t(x)$ and $n_t(y)$ for any $t$ as before,
  such that $\xi_t=n_t(x)\delta_x+n_t(y)\delta_y$
  satisfies~(\ref{eq:limit1}), which recasts into the following system
  of coupled ordinary differential equations:
\end{description}
\begin{equation}
  \label{eq:dimorph}
  \begin{aligned}
    \frac{d}{dt}n_t(x) & \!=\!n_t(x)\big(b(x,V(0)n_t(x)\!+\!V(x\!-\!y)n_t(y))
    \!-\!d(x,U(0)n_t(x)\!+\!U(x\!-\!y)n_t(y))\big) \\
    \frac{d}{dt}n_t(y) & \!=\!n_t(y)\big(b(y,V(0)n_t(y)\!+\!V(y\!-\!x)n_t(x))
    \!-\!d(y,U(0)n_t(y)\!+\!U(y\!-\!x)n_t(x))\big).
  \end{aligned}
\end{equation}

\subsection{Large-population limit with accelerated births and
deaths} \label{sec:accel}

We consider here an alternative limit of a large population,
combined with accelerated birth and death. This may  be useful to
investigate the qualitative differences of evolutionary dynamics
across  populations with allometric demographies (larger
populations made up of smaller individuals who reproduce and die
faster).

Here, we assume for simplicity that ${\cal X}= \rit^d$.  Let us
denote by $M_F$ the space $M_F(\rit^d)$.  We consider the
acceleration of birth and death processes at a rate proportional
to $K^{\eta}$ while preserving the demographic balance. That is,
the birth and death rates scale with system size according to\\

{\bf Assumption (H2):}

\begin{equation*}
  b_K(x,\zeta)=K^{\eta}r(x)+b(x,\zeta),\quad
  d_K(x,\zeta)=K^{\eta}r(x)+d(x,\zeta).
\end{equation*}
The allometric effect (smaller individuals reproduce and die
faster) is parameterized by the function $r$, positive and bounded
over ${\cal X}$, and the constant $\eta$. As in
Section~\ref{sec:limit1}, the interaction kernels $V$ and $U$ are
renormalized by $K$.  Using similar arguments as in
Section~\ref{sec:limit1}, the process $X^K={1\over K}\nu^K$ is now
a Markov process with generator
\begin{align*}
  L^K\phi(\nu)&=K\int_{\rit^d} (K^{\eta}
  r(x)+b(x,V*\nu(x)))(1-\mu_K(x))(\phi(\nu+{1\over
    K}\delta_x)-\phi(\nu)) \nu(dx) 
  \\ & +K\int_{\rit^d} (K^{\eta}
  r(x)+b(x,V*\nu(x)))\mu_K(x)\int_{\rit^d}(\phi(\nu+{1\over
    K}\delta_{z})-\phi(\nu))M_K(x,z)dz \nu(dx) 
  \\ & +K\int_{\rit^d}
  (K^{\eta} r(x)+d(x,U*\nu(x)))(\phi(\nu-{1\over
    K}\delta_x)-\phi(\nu)) \nu(dx).
\end{align*}
As before, for any measurable functions $\phi$ on $M_F$ such that
$|\phi(\nu)|+|L^K\phi(\nu)|\leq C(1+\langle\nu,1\rangle^3)$, the process
\begin{equation}
  \label{eq:mart-gal}
  \phi(X^K_t)-\phi(X^K_0)-\int_0^tL^K\phi(X^K_s)ds
\end{equation}
is a martingale. In particular, for each measurable bounded
function $f$, we obtain
\begin{align}
  & M^{K,f}_t= \langle X^K_t,f\rangle -
  \langle X^K_0,f\rangle \notag \\ &
  -\int_0^t\int_{\rit^d}(b(x,V* X^K_s(x))-d(x,U*X^K_s(x)))f(x)
  X^K_s(dx)ds \label{mart1} \\ & -\int_0^t\int_{\rit^d}
  \mu_K(x)(K^{\eta}r(x)+b(x,V*X^K_s(x)))
  \bigg(\int_{\rit^d}f(z)M_K(x,z)dz-f(x)\bigg)X^K_s(dx)ds, \notag
\end{align}
is a square integrable martingale with quadratic variation
\begin{multline}
  \label{quadra1}
  \langle M^{K,f}\rangle_t={1\over K}\bigg\{
  \int_0^t\int_{\rit^d}(2K^{\eta}r(x)+ b(x,V*X^K_s(x))+d(x,U*
  X^K_s(x)))f^2(x)X^K_s(dx)ds \\
  +\int_0^t\int_{\rit^d}\mu_K(x)(K^{\eta}r(x)+b(x,V*X^K_s(x)))
  \bigg(\int_{\rit^d}f^2(z)M_K(x,z)dz- f^2(x)\bigg)X^K_s(dx)ds
  \bigg\}.
\end{multline}

Two interesting cases will be considered hereafter, in which the
variance effect $\:\mu_K M_K\:$ is of order $1/K^{\eta}$. That
will ensure the deterministic part in~(\ref{mart1}) to converge.
In the large-population renormalization
(Section~\ref{sec:limit1}), the quadratic variation of the
martingale part was of the order of $1/K$. Here, it is of the
order of $K^{\eta}\times 1/K$. This quadratic variation will thus
stay finite provided that $\eta\in(0,1]$, in which case tractable
limits will result.  Moreover, this limit will be zero if $\eta<1$
and nonzero if $\eta=1$, which will lead to deterministic or
random limit models.

\subsubsection{Accelerated mutation and small mutation steps}
\label{sec:limit2}

We consider here that the mutation rate is fixed, so that
mutations are accelerated as a consequence of accelerating birth.
We assume\\

{\bf Assumptions (H3):}
\begin{description}
\item[\textmd{(1)}] $\mu_K=\mu$.
\item[\textmd{(2)}] The mutation step density $M_K(x,z)$ is the
  density of a random variable with mean $x$, variance-covariance
  matrix $\Sigma(x)/K^{\eta}$ (where
  $\Sigma(x)=(\Sigma_{ij}(x))_{1\leq i,j\leq d}$) and with third
  moment of order $1/K^{\eta+\varepsilon}$ uniformly in $x$
  ($\varepsilon>0$). (Thus, as $K$ goes to infinity, mutant traits
  become more concentrated around their progenitors').
\item[\textmd{(3)}] $\sqrt{\Sigma}$ denoting the symmetrical square
  root matrix of $\Sigma$, the function $\sqrt{\Sigma r\mu}$ is
  Lipschitz continuous.
\end{description}

The main example is when the mutation step density is taken as the
density of a vector of independent Gaussian variables with mean $x$
and variance $\sigma^2(x)/K^{\eta}$:
\begin{equation}
  \label{eq:gauss}
  M_K(x,z)=\left(\frac{K^{\eta}}{2\pi\sigma^2(x)}\right)^{d/2}
  \exp[-K^{\eta}|z-x|^2/2\sigma^2(x)]
\end{equation}
where $\sigma^2(x)$ is positive and bounded over $\rit^d$.

Then the convergence results of this section can be stated as follows.
\begin{thm}
  \label{readif}
  \begin{description}
  \item[\textmd{(1)}] Assume~(H), (H1), (H2), (H3) and $0<\eta<1$.
    Assume also that the initial conditions $X^K_0$ converge in law
    and for the weak topology on $M_F$ as $K$ increases, to a finite
    deterministic measure $\xi_0$, and that
    \begin{equation}
      \label{eq:X3}
      \sup_K E(\langle X^K_0,1\rangle^3)<+\infty.
    \end{equation}
    Then, for each $T>0$, the sequence of processes $(X^K)$ belonging
    to $\dit([0,T],M_F)$ converges (in law) to the unique
    deterministic function $(\xi_t)_{t\geq 0} \in C([0,T],M_F)$
    satisfying: for each function $f\in C^2_b(\rit^d)$,
    \begin{align}
      \label{readif1}
      \langle\xi_t,f\rangle &
      =\langle \xi_0,f\rangle+ \int_0^t\int_{\rit^d}
      (b(x,V*\xi_s(x))-d(x,U*\xi_s(x)))f(x)\xi_s(dx)ds \notag \\ &
      +\int_0^t\int_{\rit^d}{1\over 2}\mu(x)r(x)
      \sum_{1\leq i,j\leq d}\Sigma_{ij}(x)\partial^2_{ij}
      f(x)\xi_s(dx)ds,
    \end{align}
    where $\partial^2_{ij}f$ denotes the second-order partial
    derivative of $f$ with respect to $x_i$ and $x_j$
    ($x=(x_1,\ldots,x_d)$).
  \item[\textmd{(2)}] Assume moreover that there exists $c>0$ such
    that $r(x)\mu(x)s^*\Sigma(x)s\geq c||s||^2$ for any $x$ and $s$ in
    $\rit^d$. Then for each $t>0$, the measure $\xi_t$ has a density
    with respect to Lebesgue measure.
  \end{description}
\end{thm}

\begin{rema}
  In case (2), Eq.~(\ref {readif1}) may be written as
  \begin{equation}
    \label{readifde}
    \partial_t
    \xi_t(x)=\bigg(b(x,V*\xi_t(x))-d(x,U*\xi_t(x))\bigg)\xi_t(x) +
    {1\over 2}\sum_{1\leq i,j\leq d}
    \partial^2_{ij}(r\mu\Sigma_{ij}\xi_t)(x).
  \end{equation}
  Observe that, for the example~(\ref{eq:gauss}), this equation writes
  \begin{equation}
    \label{eq:reacdiff}
    \partial_t\xi_t(x)=\bigg(b(x,V*\xi_t(x))-d(x,U*\xi_t(x))\bigg)\xi_t(x)
    +{1\over 2}\Delta(\sigma^2r\mu\xi_t)(x).
  \end{equation}
  Therefore, Eq.~(\ref{eq:reacdiff}) generalizes the Fisher
  reaction-diffusion equation known from classical population genetics
  (see e.g.~\cite{Bu00}).
\end{rema}

\begin{thm}
  \label{readifstoch}
  Assume~(H), (H1), (H2), (H3) and $\eta=1$. Assume also that the
  initial conditions $X^K_0$ converge in law and for the weak topology
  on $M_F({\cal X})$ as $K$ increases, to a finite (possibly random)
  measure $X_0$, and that
  $\: \sup_K E(\langle X^K_0,1\rangle^3)<+\infty.$

  Then, for each $T>0$, the sequence of processes $(X^K)$ converges in
  law in $\dit([0,T],M_F)$ to the unique (in law) continuous
  superprocess $X \in C([0,T],M_F)$, defined by the following
  conditions:
  \begin{equation}
    \label{eca}
    \sup_{t \in [0,T]}E \left( \langle X_t,1\rangle^3 \right) <\infty,
  \end{equation}
  and for any $f \in C^2_b(\mathbb{R}^d)$,
  \begin{align}
    \label{dfmf}
    \bar M^f_t&=\langle X_t,f\rangle - \langle X_0,f\rangle - \frac{1}{2}
    \int_0^t  \int_{\rit^d}  \mu(x)r(x)
    \sum_{1\leq i,j\leq d}\Sigma_{ij}(x)\partial^2_{ij}
    f(x) X_s(dx)ds \notag
    \\ &- \int_0^t \int_{\rit^d} f(x)\left(
      b(x,V*X_s(x))-d(x,U*X_s(x))\right)X_s(dx)ds
  \end{align}
  is a continuous martingale with quadratic variation
  \begin{equation}
    \label{qvmf}
    \langle\bar M^f\rangle_t= 2\int_0^t \int_{\rit^d}r(x) f^2(x) X_s(dx)ds.
  \end{equation}
\end{thm}

\begin{rema}
  \begin{description}
  \item[\textmd{(1)}] The limiting measure-valued process $X$ appears
    as a generalization of the one proposed by Etheridge~\cite{Et04}
    to model spatially structured populations.
  \item[\textmd{(2)}] The conditions characterizing the process $X$
    above can be formally rewritten as equation
    \begin{equation*}
      \partial_t X_t(x) =\bigg(b(x,V*X_t(x))-d(x,U*X_t(x))\bigg)X_t(x)
      +{1\over 2}\sum_{1\leq i,j\leq d}
      \partial^2_{ij}(r\mu\Sigma_{ij}X_t)(x)+ \dot{M}_t
    \end{equation*}
    where $\dot{M}_t$ is a random fluctuation term, which reflects the
    demographic stochasticity of this fast birth-and-death process,
    that is, faster than the accelerated birth-and-death process which
    led to the deterministic reaction-diffusion
    approximation~(\ref{eq:reacdiff}).
  \item[\textmd{(3)}] As developed in Step~1 of the proof of
    Theorem~\ref{readifstoch} below, a Girsanov's theorem relates the law of
    $X_t$ and the one of a standard super-Brownian motion, which leads to
    conjecture that a density for $X_t$ exists only when $d=1$, as
    for the
    super-Brownian motion.
  \end{description}
\end{rema}

These two theorems are illustrated by the simulations of
Figs.~\ref{fig:IBM-2}~(a), (c) and~(d).

\paragraph{Proof of Theorem~\ref{readif}} We divide the proof in several
steps. Let us fix $T>0$.\bigskip

{\bf Step 1} \phantom{9} Let us first show the uniqueness for a
solution of the equation~(\ref{readif1}).

To this aim, we define the evolution equation associated
with~(\ref{readif1}). It is easy to prove that if $\xi$ is a
solution of~(\ref{readif1}) satisfying $\sup_{t\in [0,T]}
\langle\xi_t,1 \rangle <\infty$, then for each test function
$\psi_t(x)=\psi(t,x) \in C^{1,2}_b(\rit_+\times \rit^d)$, one has
\begin{align}
  \label{readif2}
  \langle\xi_t,\psi_t\rangle&=\langle \xi_0,\psi_0\rangle+
  \int_0^t\int_{\rit^d}
  (b(x,V*\xi_s(x))-d(x,U*\xi_s(x)))\psi(s,x)\xi_s(dx)ds \notag \\ &
  +\int_0^t\int_{\rit^d}(\partial_s \psi(s,x)+{1\over
    2}r(x)\mu(x)\sum_{i,j}\Sigma_{ij}(x)\partial^2_{ij}
    \psi_s(x))\xi_s(dx)ds.
\end{align}
Now, since the function $\sqrt{\Sigma r \mu}$ is Lipschitz
continuous, we may define the transition semigroup $(P_t)$ whith
infinitesimal generator $f\mapsto{1\over 2} r
\mu\sum_{i,j}\Sigma_{ij}\partial^2_{ij}f$.  Then, for each
function $f\in C^{2}_b(\rit^d)$ and fixed $t>0$, to choose
$\psi(s,x)=P_{t-s}f(x)$ yields
\begin{equation}
  \label{readif3} \langle
  \xi_t,f\rangle=\langle \xi_0,P_tf\rangle+ \int_0^t\int_{\rit^d}
  (b(x,V*\xi_s(x))-d(x,U*\xi_s(x)))P_{t-s}f(x)\xi_s(dx)ds,
\end{equation}
since $\partial_s \psi(s,x)+{1\over
  2}r(x)\mu(x)\sum_{i,j}\Sigma_{ij}(x)\partial^2_{ij} \psi_s(x)=0$ for
this choice.

We now prove the uniqueness of a solution of~(\ref{readif3}).

Let us consider two solutions $(\xi_t)_{t\geq 0}$ and $(\bar
\xi_t)_{t\geq 0}$ of~(\ref{readif3}) satisfying $\sup_{t\in [0,T]}
\left< \xi_t+\bar \xi_t,1 \right> =A_T <+\infty$. We consider the
variation norm defined for $\mu_1$ and $\mu_2$ in $M_F$ by
\begin{equation}
  || \mu_1 - \mu_2 || = \sup_{f \in
    L^\infty({\rit^d}), \; | | f ||_\infty \leq 1}
  |\left<\mu_1-\mu_2,f\right>|.
\end{equation}
Then, we consider some bounded and measurable function $f$ defined
on ${\cal X}$ such that $|| f ||_\infty \leq 1$ and obtain
\begin{align}
  \label{lip}
  | \left<\xi_t-\bar \xi_t,f\right> |&\leq
  \intot \left| \int_{\rit^d} [\xi_s(dx) - \bar \xi_s(dx)]
    \left(b(x,V*\xi_s(x))-d(x,U*\xi_s(x))\right)P_{t-s}f(x)\right| ds
    \notag \\ & + \intot \left|\int_{\rit^d} \bar
    \xi_s(dx)(b(x,V*\xi_s(x))-b(x,V*\bar{\xi}_s(x)))P_{t-s}f(x)\right|
    ds \notag \\ & + \intot \left|\int_{\rit^d} \bar
    \xi_s(dx)(d(x,U*\xi_s(x))-d(x,U*\bar{\xi}_s(x)))P_{t-s}f(x)\right|
    ds.
\end{align}

Since $|| f ||_\infty \leq 1$, then $|| P_{t-s}f ||_\infty \leq 1$
and for all $x\in \rit^d$,
\begin{equation*}
  \left|(b(x,V*\xi_s(x))-d(x,U*\xi_s(x)))P_{t-s}f(x)\right|\leq
  \bar{b}+\bar{d}(1+\bar{U}A_T).
\end{equation*}
Moreover, $b$ and $d$ are Lipschitz continuous in their second
variable with respective constants $K_b$ and $K_d$. Thus we obtain
from~(\ref{lip}) that
\begin{equation}
  | \left<\xi_t-\bar \xi_t,f\right> | \leq
  \left[\bar{b}+\bar{d}(1+\bar{U}A_T) +K_b A_T\bar{V}+K_d A_T\bar{U}
  \right] \intot  || \xi_s - \bar\xi_s ||ds.
\end{equation}
Taking the supremum over all functions $f$ such that $| | f
||_\infty \leq 1$, and using the Gronwall Lemma, we finally deduce
that for all $t\leq T$, $|| \xi_t - \bar\xi_t| | =0$. Uniqueness
holds.\bigskip

{\bf Step 2} \phantom{9} Next, we would like to obtain some moment
estimates.  First, we check that for all $T<\infty$,
\begin{equation}
  \label{me1}
  \sup_K \sup_{t \in [0,T]}  E \big(\langle X^K_t,1\rangle^3
  \big)
  <\infty.
\end{equation}
To this end, we use~(\ref{eq:mart-gal}) with $\phi(\nu)=\langle
\nu,1\rangle^3$. (To be completely rigorous, one should first use
$\phi(\nu)=\langle \nu,1\rangle^3 \land A$, make $A$ tend to
infinity). Taking expectation, we obtain that for all $t\geq 0$,
all $K$,
\begin{align*}
  &E \left(\langle X^K_t,1\rangle^3 \right) =E \left(\langle
    X^K_0,1\rangle^3 \right) \\ &+ \intot  E \bigg(\int_{\rit^d}
   \bigg([K^{\eta+1}r(x) +K b(x,V* X^K_s(x))] \left\{[\langle
      X^K_s,1\rangle+ \frac{1}{K}]^3-\langle X^K_s,1\rangle^3 \right\}
  \\ &\left\{K^{\eta+1}r(x) + K d(x,U* X^K_s(x)) \right\}
    \left\{[\langle X^K_s,1\rangle- \frac{1}{K}]^3-\langle
      X^K_s,1\rangle^3 \right\}\bigg) X^K_s(dx) \bigg)ds.
\end{align*}
Neglecting the non-positive death term involving $d$, we get
\begin{align*}
  E & \left(\langle X^K_t,1\rangle^3 \right) \leq E \left(\langle
    X^K_0,1\rangle^3 \right) \\ &+ \intot E \bigg(\int_{\rit^d}\bigg(
  K^{\eta+1}r(x) \left\{[\langle X^K_s,1\rangle+ \frac{1}{K}]^3 +
    [\langle X^K_s,1\rangle- \frac{1}{K}]^3 -2\langle X^K_s,1\rangle^3
  \right\} \\ & + Kb(x,V* X^K_s(x)) \left\{[\langle X^K_s,1\rangle+
    \frac{1}{K}]^3-\langle X^K_s,1\rangle^3 \right\}\bigg)X^K_s(dx)
  \bigg)ds.
\end{align*}
But for all $x\geq 0$, all $\e \in (0,1]$, $(x+\e)^3 -x^3 \leq 6
\e (1+x^2)$ and $| (x+\e)^3+(x-\e)^3-2x^3 | = 6 \e^2 x$. We
finally obtain
\begin{equation*}
  E \left(\langle X^K_t,1\rangle^3 \right) \leq E \left(\langle
      X^K_0,1\rangle^3 \right)+ C \intot  E\left(\langle
      X^K_s,1\rangle +\langle X^K_s,1\rangle^2  +\langle X^K_s,1\rangle^3
  \right)ds.
\end{equation*}
Assumption~(\ref{eq:X3}) and
the Gronwall Lemma allows us to conclude that~(\ref{me1}) holds.\\
Next, we wish to check that
\begin{equation}
  \label{me2}
  \sup_K E \big(\sup_{t \in [0,T]} \langle X^K_t,1\rangle^2\big)
  <\infty.
\end{equation}
Applying~(\ref{mart1}) with $f \equiv 1$, we obtain
\begin{equation*}
  \langle X^K_t,1\rangle = \langle X^K_0,1\rangle+ \intot  \int_{\cal X}
  \left( b(x,V* X^K_s(x))-d(x,U* X^K_s(x)) \right) X^K_s(dx)ds +
  m^{K,1}_t.
\end{equation*}
Hence
\begin{equation*}
  \sup_{s \in [0,t] }\langle X^K_s,1\rangle^2 \leq C\bigg(\langle
    X^K_0,1\rangle^2 + \bar{b} \intot  \langle X^K_s,1\rangle^2 ds +
  \sup_{s \in [0,t] } | M^{K,1}_s |^2 \bigg).
\end{equation*}
Thanks to~(\ref{eq:X3}), the Doob inequality and the Gronwall
Lemma, there exists a constant $C_t$ not depending on $K$ such
that
\begin{equation*}
  E \big(\sup_{s \in [0,t]} \langle X^K_s,1\rangle^2 \big) \leq
  C_t \left(1+ E \left( \langle M^{K,1} \rangle_t \right) \right).
\end{equation*}
Using now~(\ref{quadra1}), we obtain, for some other constant
$C_t$ not depending on $K$,
\begin{equation*}
  E \left( \langle M^{K,1} \rangle_t \right) \leq C \intot
  \big(E\left( \langle X^K_s,1\rangle +
    \langle X^K_s,1\rangle^2 \right)\big)ds \leq C_t
\end{equation*}
thanks to~(\ref{me1}). This concludes the proof
of~(\ref{me2}).\bigskip

{\bf Step 3} \phantom{9} We first endow $M_F$ with the vague
topology, the extension to the weak topology being handled in Step
6 below. To show the tightness of the sequence of laws
$Q^K=\loi(X^K)$ in ${\cal
  P}(\dit([0,T],M_F))$, it suffices, following Roelly~\cite{Ro86}, to
show that for any continuous bounded function $f$ on $\rit^d$, the
sequence of laws of the processes $\langle X^K,f\rangle$ is tight
in $\dit([0,T], \rit)$. To this end, we use the Aldous
criterion~\cite{Al78} and the Rebolledo criterion
(see~\cite{JM86}). We have to show that
\begin{equation}
  \label{tight} \sup_K
  E\big(\sup_{t\in [0,T]} | \langle X^K_s,f\rangle | \big)
  <\infty,
\end{equation}
and the tightness respectively of the laws of the predictable
quadratic variation of the martingale part and
of the drift part of the semimartingales $\langle X^K,f\rangle$. \\
Since $f$ is bounded,~(\ref{tight}) is a consequence
of~(\ref{me2}): let us thus consider a couple $(S,S')$ of stopping
times satisfying a.s.~$0 \leq S \leq S' \leq S+\delta\leq T$.
Using~(\ref{quadra1}) and (\ref{me2}), we get for constants $C,
C'$
\begin{equation*}
  E\left(\langle M^{K,f}\rangle_{S'}- \langle
    M^{K,f}\rangle_S\right) \leq C E\left( \int_S^{S+\delta} \left(
      \langle X^K_s,1\rangle+ \langle X^K_s,1\rangle^2
    \right)ds\right)\leq C' \delta.
\end{equation*}
In a similar way,  the expectation of the finite variation part of
$\langle X^K_{S'},f\rangle - \langle X^K_S,f\rangle$ is bounded by
$C' \delta$.

Hence,  the sequence $Q^K=\loi(X^K)$ is tight.\bigskip

{\bf Step 4} \phantom{9} Let us now denote by $Q$ the limiting law
of a subsequence of $Q^K$. We still denote this subsequence by
$Q^K$. Let $X=(X_t)_{t\geq 0}$ a process with law $Q$. We remark
that by construction, almost surely,
\begin{equation*}
  \sup_{t\in [0,T]}\ \sup_{f \in L^\infty(\rit^d), || f
    ||_\infty\leq 1} | \langle X^K_t,f \rangle -
  \langle X^K_{t^-},f \rangle | \leq 1/K.
\end{equation*}
This implies that the process $X$ is a.s.~strongly continuous.
\bigskip

{\bf Step 5} \phantom{9} The time $T>0$ is fixed. Let us now check
that almost surely, the process $X$ is the unique solution
of~(\ref{readif1}).  Thanks to~(\ref{me2}), it satisfies
$\sup_{t\in
  [0,T]} \langle X_t,1\rangle <+\infty$ a.s., for each $T$. We fix now
a function $f\in C^3_b(\rit^d)$ (the extension of~(\ref{readif1})
to
any function $f$ in $C^2_b$ is not hard) and some $t\leq T$.\\
For $\nu \in C([0,T],M_F)$, denote by
\begin{align}
  \Psi^1_t(\nu)&= \langle\nu_t,f\rangle - \langle\nu_0,f\rangle -
  \intot \int_{\rit^d}
  (b(x,V*\nu_s(x))-d(x,U*\nu_s(x)))f(x)\nu_s(dx)ds,\notag \\
  \Psi^2_t(\nu)&=-\intot \int_{\rit^d}  {1\over
    2}\mu(x)r(x)\sum_{i,j}\Sigma_{ij}(x)\partial^2_{ij}f(x)\nu_s(dx) ds.
\end{align}
We have to show that
\begin{equation}
  \label{wwhtp}
  E_Q \left( |\Psi^1_t(X)+\Psi^2_t(X) | \right)=0.
\end{equation}

By~(\ref{mart1}), we know that for each $K$,
\begin{equation*}
  M^{K,f}_t=\Psi^1_t(X^K)+\Psi^{2,K}_t(X^K),
\end{equation*}
where
\begin{multline}
  \Psi^{2,K}_t(X^K)=-\int_0^t\int_{\rit^d}
  \mu(x)(K^{\eta}r(x)+b(x,V*X^K_s(x)))\\ \bigg(\int_{\rit^d}
  f(z)M_K(x,z)dz-f(x)\bigg)X^K_s(dx)ds.
\end{multline}
Moreover,~(\ref{me2}) implies that for each $K$,
\begin{equation}
  \label{cqmke}
  E \left( | M^{K,f}_t |^2 \right) = E \left( \langle
      M^{K,f}\rangle_t \right)\leq \frac{C_{f}K^{\eta}}{K}E\left(
    \intot  \left\{\langle X^K_s,1\rangle
      +\langle X^K_s,1\rangle^2\right\}ds \right) \leq
  \frac{C_{f,T}K^{\eta}}{K},
\end{equation}
which goes to $0$ as $K$ tends to infinity, since $0<\alpha<1$.
Therefore,
\begin{equation*}
  \lim_K E(|\Psi^1_t(X^K)+\Psi^{2,K}_t(X^K)|)=0.
\end{equation*}

Since $X$ is a.s.~strongly continuous, since $f\in C^3_b(\rit^d)$
and thanks to the continuity of the parameters, the functions
$\Psi^1_t$ and $\Psi^2_t$ are a.s.~continuous at $X$. Furthermore,
for any $ \nu \in \dit([0,T],M_F)$,
\begin{equation}
  |\Psi^1_t(\nu)+\Psi^2_t(\nu)| \leq C_{f,T} \sup_{s\in [0,T]}
  \left(1+\langle \nu_s,1\rangle^2 \right).
\end{equation}
Hence using~(\ref{me1}), we see that the sequence
$(\Psi^1_t(X^K)+\Psi^2_t(X^K))_K$ is uniformly integrable, and
thus
\begin{eqnarray}
  \label{cqv1}
  \lim_K E\left(|\Psi^1_t(X^K)+\Psi^2_t(X^K)|\right)
  =E\left(|\Psi^1_t(X)+\Psi^2_t(X)|\right).
\end{eqnarray}

We have now to deal with $\Psi^{2,K}_t(X^K)-\Psi^2_t(X^K)$. The
convergence of this term is due to the fact that the measure
$M_K(x,z)dz$ has mean $x$, variance $\Sigma(x)/K^{\eta}$, and
third moment bounded by $C/K^{\eta+\varepsilon}$ ($\varepsilon>0$)
uniformly in $x$.  Indeed, if $Hf(x)$ denotes the Hessian matrix
of $f$ at $x$,
\begin{align}
  & \int_{\rit^d} f(z) M_K(x,z)dz \notag \\ &=\int_{\rit^d}
  \left(f(x)+(z-x)\cdot\nabla f(x) +\frac{1}{2} (z-x)^*
  Hf(x)(z-x)+O((z-x)^3)\right) M_K(x,z)dz \notag \\ &= f(x)+{1\over
    2}\sum_{i,j}\frac{\Sigma_{ij}(x)}{K^{\eta}}\partial^2_{ij}
  f(x)+{\it o}({1\over K^{\eta}}). \label{eq:M_K}
\end{align}
where $K^{\eta}{\it o}({1\over K^{\eta}})$ tends to $0$ uniformly
in $x$ (since $f$ is in $C^3_b$), as $K$ tends to infinity.  Then,
\begin{multline*}
  \Psi^{2,K}_t(X^K)=-\int_0^t\int_{\rit^d}
  \mu(x)(K^{\eta}r(x)+b(x,V*X^K_s(x)))\times \\ \times\bigg({1\over 2}
  \sum_{i,j}\frac{\Sigma_{ij}(x)}{K^{\eta}}\partial^2_{ij}
  f(x)+{\it o}({1\over K^{\eta}})\bigg)X^K_s(dx)ds,
\end{multline*}
and
\begin{equation*}
  |\Psi^{2,K}_t(X^K)-\Psi^2_t(X^K)| \leq C_f
  <X^K_s,1>\bigg({1\over K^{\eta}}+K^{\eta}{\it o}({1\over
    K^{\eta}})\bigg).
\end{equation*}
Using~(\ref{me2}), we conclude the proof of~(\ref{wwhtp}).\bigskip

{\bf Step 6} \phantom{9} The previous steps imply that $(X^K)_K$
converges to $\xi$ in $\dit([0,T],M_F)$, where $M_F$ is endowed with
the vague topology. To extend the result to the case where $M_F$ is
endowed with the weak topology, we use a criterion proved in
M\'el\'eard and Roelly~\cite{MR93}: since the limiting process is
continuous, it suffices to prove that the sequence $(\langle
X^K,1\rangle)$ converges to $\langle\xi,1\rangle$ in law, in
$\dit([0,T],\rit)$. One may of course apply Step 5 with $f\equiv 1$,
which concludes the proof.  \bigskip

2) Let us now assume the non-degeneracy property
$r(x)\mu(x)s^*\Sigma(x)s\geq c\|s\|^2>0$ for each $x\in \rit^d, s
\in\rit^d$. That implies that for each time $t>0$, the transition
semigroup $P_t(x,dy)$ introduced in Step 1 of this proof has for
each $x$ a density function $p_t(x,y)$ with respect to the
Lebesgue measure. Then if we come back to the evolution
equation~(\ref{readif3}), we can write
\begin{multline*}
  \int_{\rit^d} f(x) \xi_t(dx) =\int_{\rit^d}
  \left(\int_{\rit^d}f(y)p_t(x,y)dy\right) \xi_0(dx)\\ +
  \int_0^t\int_{\rit^d}
  (b(x,V*\xi_s(x))-d(x,U*\xi_s(x)))\bigg(\int_{\rit^d}
  f(y)p_{t-s}(x,y)dy\bigg)\xi_s(dx)ds.
\end{multline*}

Using the fact that the parameters are bounded, that $\
\sup_{t\leq
  T}\langle\xi_t,1\rangle<+\infty\ $ and that $\ f\ $ is bounded, we can
apply Fubini's theorem and deduce that
\begin{equation*}
  \int_{\rit^d} f(x) \xi_t(dx)=\int_{\rit^d}H_t(y)f(y)dy
\end{equation*}
with $H\in L^{\infty}([0,T],L^1(\rit^d))$, which implies that
$\xi_t$ has a density with respect to the Lebesgue measure for
each time $t\leq T$.

Equation~(\ref{readifde}) is then the dual form of
(\ref{readif1}).\hfill$\Box$
\bigskip

\paragraph{Proof of Theorem~\ref{readifstoch}}
We will use a similar method as the one of the previous theorem.
Steps~2, 3, 4 and~6 of this proof can be achieved exactly in the
same way. Therefore, we only have to prove the uniqueness (in law)
of the solution to the martingale
problem~(\ref{eca})--(\ref{qvmf}) (Step~1), and that any
accumulation point of the sequence of laws of $X^K$ is solution
to~(\ref{eca})--(\ref{qvmf}) (Step~5).\bigskip

{\bf Step 1} \phantom{9} This uniqueness result is well-known for
the super-Brownian process (defined by a similar martingale
problem, but with $b=d=0$, $r=\mu=1$ and $\Sigma=\mbox{Id}$, cf.
\cite{Ro86}). Following ~\cite{Et04}, we may use the version of
Dawson's Girsanov transform obtained in Evans and
Perkins~\cite{EP94} (Theorem~2.3), to deduce the uniqueness in our
situation, provided the condition
\begin{equation*}
  E\left(\intot \intrd
    [b(x,V*X_s(x))-d(x,U*X_s(x))]^2X_s(dx) ds\right)<+\infty
\end{equation*}
is satisfied. This is easily obtained from the assumption that
$\sup_{t\in [0,T]}E[\langle X_t,1\rangle^3]<\infty$ since the
coefficients are bounded.\bigskip

{\bf Step 5} \phantom{9} Let us identify the limit. Let us call
$Q^K=\loi(X^K)$ and denote by $Q$ a limiting value of the tight
sequence $Q^K$, and by $X=(X_t)_{t\geq 0}$ a process with law $Q$.
Because of Step~4, $X$ belongs a.s.~to $C([0,T],M_F)$. We have to
show that $X$ satisfies the conditions~(\ref{eca}), (\ref{dfmf})
and~(\ref{qvmf}). First note that~(\ref{eca}) is straightforward
from~(\ref{me2}). Then, we show that for any function $f$ in
$C^3_b(\mathbb{R}^d)$, the process $\bar M^f_t$ defined
by~(\ref{dfmf}) is a martingale (the extension to every function
in $C^2_b$ is not hard). We consider $0\leq s_1\leq...\leq
s_n<s<t$, some continuous bounded maps $\phi_1,...\phi_n$ on
$M_F$, and our aim is to prove that, if the function $\Psi$ from
$\mathbb{D}([0,T],M_F)$ into $\mathbb{R}$ is defined by
\begin{align}
  &\Psi(\nu) = \phi_1(\nu_{s_1})...\phi_n(\nu_{s_n}) \Big\{
  \langle \nu_t,f\rangle -\langle \nu_s,f\rangle
  \notag \\ & -\! \int_s^t\! \intrd\! \bigg({1\over 2}\mu(x)
  r(x)\sum_{i,j}\Sigma_{ij}\partial^2_{ij} f(x)+ f(x)
  \left[b(x,V*\nu_u(x))-d(x,U*\nu_u(x)) \right]\bigg)\nu_u(dx)du
  \Big\},
\end{align}
then
\begin{equation}
  \label{cqfd44}
  E\left( \Psi(X) \right)=0.
\end{equation}
It follows from~(\ref{mart1}) that
\begin{eqnarray}
  \label{262626}
  0=E \left(\phi_1(X^K_{s_1})...\phi_n(X^K_{s_n})\left\{ M^{K,f}_t -
    M^{K,f}_s \right\}\right)=E \left( \Psi(X^K) \right) - A_K,
\end{eqnarray}
where $A_K$ is defined by
\begin{multline*}
  A_K=E \Big( \phi_1(X^K_{s_1})...\phi_n(X^K_{s_n})\int_s^t \intrd
    \mu(x) \Big\{ b(x,V*X^K_u(x))\Big[\intrd  (f(z)- f(x))M_K(x,z)dz\Big]\\
    +r(x) K\Big[ \intrd  (f(z) - f(x)
    - \sum_{i,j}{\Sigma_{ij}(x)\over
    2K}\partial^2_{ij}f(x))M_K(x,z)dz\Big]
  \Big\}X^K_u(dx)du\Big).
\end{multline*}
It turns out from~(\ref{eq:M_K}) that $A_K$ tends to zero  as $K$
grows to infinity, and using~(\ref{me2}), that the sequence
$(|\Psi(X^K)|)_K$ is uniformly integrable, so
\begin{equation}
  \label{333}
  \lim_K E\left(|\Psi (X^K)|\right) = E_Q\left(|\Psi(X)| \right).
\end{equation}
Collecting the previous results allows us to conclude
that~(\ref{cqfd44}) holds, and thus $\bar M^f$ is a martingale.\\
We finally have to show that the bracket of $\bar M^f$ is given
by~(\ref{qvmf}). To this end, we first check that
\begin{align}
  \label{lfalqoc}
  \bar N^{f}_t &= \langle X_t,f \rangle^2 - \langle X_0,f \rangle^2
  -\intot  \intrd 2r(x)f^2(x) X_s(dx)ds \notag \\
  &- 2\intot \langle X_s,f \rangle  \intrd  f(x)
  \left[b(x,V*X_s(x))-d(x,U*X_s(x)) \right]X_s(dx) ds \notag \\
  &- \intot   \langle X_s,f \rangle  \intrd
  \mu(x)r(x)\sum_{i,j}\Sigma_{ij}(x)\partial^2_{ij}f(x) X_s(dx) ds
\end{align}
is a martingale. This can be done exactly as for $\bar M^f_t$,
 using the semimartingale
decomposition of $\langle X^K_t,f \rangle^2$, given by
(\ref{eq:mart-gal}) with $\phi(\nu)=\langle\nu,f\rangle^2$.
  In another hand, It\^o's formula implies that
\begin{multline*}
  \langle X_t,f \rangle^2 - \langle X_0,f \rangle^2 - \langle \bar
  M^f\rangle_t
  - \intot ds 2 \langle X_s,f \rangle \intrd X_s(dx)\frac{1}{2}
  r(x)\mu(x)\sum_{i,j}\Sigma_{ij}(x)\partial^2_{ij}f(x) \\
  - \intot ds 2 \langle X_s,f \rangle  \intrd X_s(dx) f(x)
  \big[b(x,V*X_s(x))-d(x,U*X_s(x)) \big]
\end{multline*}
is a martingale. Comparing this formula with~(\ref{lfalqoc}), we
obtain (\ref{qvmf}).\hfill$\Box$
\bigskip

\subsubsection{Rare mutations} \label{sec:limit3}

In this case, the mutation step density $M$ is fixed and the
mutation rate is decelerated proportionally to $1/K^{\eta}$:\\

{\bf Assumption (H4):}

\begin{equation*}
  M_K=M,\quad \mu_K={\mu\over K^{\eta}}.
\end{equation*}
Thus only births without mutation are accelerated.

As in Section~\ref{sec:limit2}, we obtain deterministic or random
limits, according to the value of $\eta\in (0,1]$.
\begin{thm}
  \label{intdif}
  \begin{description}
  \item[\textmd{(1)}] Assume~(H), (H1), (H2), (H4) and $0<\eta<1$.
    Assume also that the initial conditions $X^K_0$ converge in law
    and for the weak topology on $M_F({\cal X})$ as $K$ increases, to
    a finite deterministic measure $\xi_0$, and that $\sup_K E(\langle
    X^K_0,1\rangle^3)<+\infty$.

    Then, for each $T>0$, the sequence of processes $(X^K)$ belonging
    to $\dit([0,T],M_F)$ converges (in law) to the unique
    deterministic function $(\xi_t)_{t\geq 0} \in C([0,T],M_F)$ weak
    solution of the deterministic nonlinear integro-differential
    equation:
    \begin{equation}
      \label{eq:intdif}
      \partial_t\xi_t(x)=[b(x,V*\xi_t(x))-d(x,U*\xi_t(x))]\xi_t(x)
      +\int_{\rit^d}M(y,x) \mu(y) r(y)\xi_t(y)dy.
    \end{equation}
  \item[\textmd{(2)}] Assume now $\eta=1$ and that $X^K_0$ converge in
    law to $X_0$. Then, for each $T>0$, the sequence of processes
    $(X^K)$ converges in law in $\dit([0,T],M_F)$ to the unique (in
    law) continuous superprocess $X \in C([0,T],M_F)$, defined by the
    following conditions:
    \begin{equation*}
      \sup_{t \in [0,T]}E \left( \langle X_t,1\rangle^3 \right) <\infty,
    \end{equation*}
    and for any $f \in C^2_b(\mathbb{R}^d)$,
    \begin{align*}
      \bar M^f_t & =\langle X_t,f\rangle - \langle X_0,f\rangle -
      \int_0^t  \int_{\rit^d}  \mu(x)r(x)
      \int_{\rit^d} M(x,z) f(z)dz X_s(dx)ds 
      \\ & -\int_0^t \int_{\rit^d} f(x)\left(
        b(x,V*X_s(x))-d(x,U*X_s(x))\right)X_s(dx)ds
    \end{align*}
    is a continuous martingale with quadratic variation
    \begin{equation*}
      \langle\bar M^f\rangle_t= 2\int_0^t \int_{\rit^d}r(x) f^2(x) X_s(dx)ds.
    \end{equation*}
  \end{description}
\end{thm}

In a SPDE formalism, one can write the last limit as formal
solution of the equation
\begin{equation}
  \label{eq:superproc2}
  \partial_t X_t(x)=[b(x,V*X_t(x))-d(x,U*X_t(x))]X_t(x)
  +\int_{\rit^d}M(y,x) \mu(y) r(y)X_t(dy)+\dot{M},
\end{equation}
where $\dot{M}$ is a random fluctuation term.

The proof of Theorem~\ref{intdif} is similar to proofs of
Theorems~\ref{readif} and~\ref{readifstoch} and we leave it to the
reader. Theorem~\ref{intdif}~(1) is illustrated in the simulation of
Fig.~\ref{fig:IBM-2}~(b).

\section{Rare mutation renormalization of the monomorphic process and
  adaptive dynamics}
\label{sec:AdaDyn}

In the previous section, Eqs.~(\ref{eq:intdif})
and~(\ref{eq:superproc2}) have be obtained at the population
growth time scale (ecological time scale), under an assumption of
rare mutation. Here, we are interested in the behavior of the
population process at the evolutionary time scale, when mutations
are extremely rare, as illustrated by the simulation of
Fig.~\ref{fig:IBM-1}~(d). We hence recover rigorously the
stochastic ``trait substitution sequence'' jump process of
adaptive dynamics (Metz et al.~\cite{Mal96}) when the initial
condition is monomorphic. The biological idea behind such a
scaling of the population process is that selection has sufficient
time between two mutations to eliminate all disadvantaged traits,
so that the population remains monomorphic on the evolutionary
timescale. Then the evolution proceeds by successive invasions of
mutant traits, replacing the resident trait from which the mutant
trait is born, occuring on an infinitesimal timescale with respect
to the mutation timescale. Our result emphasizes how the mutation
scaling should compare to the system size ($K$) in order to obtain
the correct time scale separation between the ``mutant-invasions''
(taking place on a short time scale) and the mutations
(evolutionary time scale).

\subsection{Statement of the result}
\label{sec:TSS}

We consider here a limit of rare mutations combined with the large
population limit of Section~\ref{sec:limit1} (Assumption~(H1) and
$b_K=b$, $d_K=d$ and $M_K=M$). We assume

{\bf Assumptions (H5):}

(i) $\mu_K(x)=u_K\mu(x)$.

(ii) For any constant $C>0$,
\begin{equation}
  \label{eq:K-u_K}
  e^{-CK} \ll u_K \ll \frac{1}{K\log K}
\end{equation}
(thus $u_K\rightarrow 0$ when $K\rightarrow +\infty$), or,
equivalently, for any $C$ and $t>0$,
\begin{equation}
  \label{eq:K-u_K-bis}
  \log K \ll \frac{t}{Ku_K} \ll e^{CK}.
\end{equation}

(iii) For any $x\in{\cal X}$, $\zeta\mapsto b(x,\zeta)$ and
$\zeta\mapsto d(x,\zeta)$ are positive functions, non-increasing and
increasing respectively, satisfying
\begin{gather}
  \forall x\in{\cal X},\ b(x,0)-d(x,0)>0, \notag \\
  \lim_{\zeta\rightarrow+\infty}\:\inf_{x\in{\cal
    X}}d(x,\zeta)=+\infty. \label{eq:d-infty}
\end{gather}

(iv) There exists a constant $\underline{U}>0$ such that
$U(h)\geq\underline{U}$ for any $h\in\mathbb{R}^d$.\\

Assumption (H5)-(i) entails the rare mutation asymptotic, and
(H5)-(ii) gives the correct scaling between the mutation
probability and the system size in order to obtain the correct
time scale separation. Observe that (H5)-(ii) implies that
$Ku_K\rightarrow 0$ when $K\rightarrow +\infty$, so that the
timescale $t/Ku_K$, which corresponds to the timescale of
mutations (the population size is proportional to $K$, and each
birth event produces a mutant with a probability proportional to
$u_K$, which gives a total mutation rate in the population
proportional to $Ku_K$) is a long timescale. Our result gives the
behavior of the population process on this long timescale.

Assumptions (H5)-(iii) and (iv) will allow to bound the population
size on the mutation timescale, and to study the behavior of the
population when it is monomorphic or dimorphic between two (rare)
mutation events. Specifically, the monotonicity properties of $b$
and $d$ in Assumption (H5)-(iii) ensures, for any $x\in{\cal X}$,
the existence of a unique non-trivial stable equilibrium
$\bar{n}(x)$ for the monomorphic logistic
equation~(\ref{eq:monomorph}) of Example~3 in
Section~\ref{sec:limit1}. Moreover, since
$b(x,V(0)u)-d(x,U(0)u)>0$ for any $u<\bar{n}(x)$ and
$b(x,V(0)u)-d(x,U(0)u)<0$ for any $u>\bar{n}(x)$, any solution
to~(\ref{eq:monomorph}) with positive initial condition converges
to $\bar{n}(x)$.

Concerning the dimorphic logistic equations~(\ref{eq:dimorph}), an
elementary linear analysis of the equilibrium $(\bar{n}(x),0)$ gives
that it is stable if $f(y,x)<0$ and unstable if $f(y,x)>0$, where the
function
\begin{equation}
  \label{eq:def-fitness}
  f(y,x)=b(y,V(y-x)\bar{n}(x))-d(y,U(y-x)\bar{n}(x))
\end{equation}
is known as the ``fitness function'' (\cite{MNG92,Mal96}), which gives
a measure of the selective advantage of a mutant individual with trait
$y$ in a monomorphic population of trait $x$ at equilibrium.
Similarly, the stability of the equilibrium $(0,\bar{n}(y))$ is
governed by the sign of $f(x,y)$.

In order to ensure that, when the invasion of a mutant trait is
possible, then this invasion will end with the extinction of the
resident trait, we will need the following additional assumption:

{\bf Assumptions (H6):}

Given any $x\in{\cal X}$, Lebesgue almost any $y\in{\cal X}$ satisfies
one of the two following conditions:

(i) either $f(y,x)<0$ (so that $(\bar{n}(x),0)$ is stable),

(ii) or $f(y,x)>0$, $f(x,y)<0$ and any solution to~(\ref{eq:dimorph})
with initial condition with positive coordinates in a given
neighborhood of $(\bar{n}(x),0)$ converges to
$(0,\bar{n}(y))$.\\

In the case of linear logistic density-dependence introduced in
Section~\ref{sec:simul} ($b(x,\zeta)=b(x)$ and
$d(x,\zeta)=d(x)+\alpha(x)\zeta$), the equilibrium monomorphic density
$\bar{n}(x)$ writes $(b(x)-d(x))/\alpha(x)U(0)$ and the
condition~(H6)-(ii) is actually equivalent to $f(y,x)>0$ and $f(x,y)<0$
(see~\cite{Ch04}).

Our convergence result writes
\begin{thm}
  \label{thm:IPS-TSS}
  Assume~(H), (H1), (H5) and~(H6). Given $x\in{\cal X}$, $\gamma>0$ and
  a sequence of $\mathbb{N}$-valued random variables
  $(\gamma_K)_{K\in\mathbb{N}}$, such that $\gamma_K/K$ is bounded in
  $\mathbb{L}^1$ and converges in law to $\gamma$, consider the
  process $(X^K_t,t\geq 0)$ of Section~\ref{sec:large-popu} generated
  by~(\ref{eq:def-LK}) with initial state
  $\frac{\gamma_K}{K}\delta_x$. Then, for any $n\geq 1$,
  $\varepsilon>0$ and $0<t_1<t_2<\ldots<t_n<\infty$, and for any
  measurable subsets $\Gamma_1,\ldots,\Gamma_n$ of ${\cal X}$,
  \begin{multline}
    \label{eq:def-cvgce}
    \lim_{K\rightarrow+\infty}P\bigl(\forall i\in\{1,\ldots,n\},\
    \exists x_i\in\Gamma_i\: :\:
    \mbox{\textnormal{Supp}}(X^K_{t_i/K u_K})=\{x_i\} \\ \mbox{and\ }|\langle
    X^K_{t_i/K u_K},\mathbf{1}\rangle-\bar{n}(x_i)|<\varepsilon\bigr)
    =P(\forall i\in\{1,\ldots,n\},\ Y_{t_i}\in\Gamma_i)
  \end{multline}
  where for any $\nu\in M_F({\cal X})$, $\mbox{\textnormal{Supp}}(\nu)$ is
  the support of $\nu$ and $(Y_t,t\geq 0)$ is a Markov jump process
  with initial state $x$ generated by
  \begin{equation}
    \label{eq:TSS-gene}
    A\varphi(x)=\int_{\rit^d}(\varphi(y)-\varphi(x))g(y,x)M(x,y)dy
  \end{equation}
  where
  \begin{equation}
    \label{eq:g}
    g(y,x)=\mu(x)b(x,V(0)\bar{n}(x))\bar{n}(x)
    \frac{[f(y,x)]_+}{b(y,V(y-x)\bar{n}(x))}
  \end{equation}
  and $[\cdot]_+$ denotes the positive part.
\end{thm}

\begin{cor}
  \label{cor:IPS-TSS}
  With the same notations and assumptions as in
  Theorem~\ref{thm:IPS-TSS}, assuming moreover that $\gamma_K/K$ is
  bounded in $\mathbb{L}^q$ for some $q>1$, the process
  $(X^K_{t/Ku_K},t\geq 0)$ converges when $K\rightarrow +\infty$, in
  the sense of the finite dimensional distributions for the topology
  on $M_F({\cal X})$ induced by the functions
  $\nu\mapsto\langle\nu,f\rangle$ with $f$ bounded and measurable on
  ${\cal X}$, to the process $(Z_t,t\geq 0)$ defined by
  \begin{equation*}
    Z_t=\left\{\begin{array}{ll}
        \gamma\delta_x & \mbox{if\ }t=0 \\
        \bar{n}(Y_t)\delta_{Y_t} & \mbox{if\ }t>0.
      \end{array}\right.
  \end{equation*}
\end{cor}

This corollary follows from the following long time moment
estimates.
\begin{lem}
  \label{lem:moment}
  Under~(H), (H1), (H5)(iii)~(\ref{eq:d-infty}) and (iv), and if $\: \sup_{K\geq 1} E(\langle
  X^K_0,1\rangle^q)<+\infty$ for some $q\geq 1$, then
  \begin{equation*}
    \sup_{K\geq 1}\:\sup_{t\geq 0}E\big(\langle
    X^K_t,\mathbf{1}\rangle^q\big)<+\infty,
  \end{equation*}
  and therefore, if $q>1$, the family of random variables $\{\langle
  X^K_t,\mathbf{1}\rangle\}_{\{K\geq 1,\:t\geq 0\}}$ is uniformly
  integrable.
\end{lem}

\paragraph{Proof of Lemma~\ref{lem:moment}}
Observe that, if we replace $b(x,V*\nu)$ by $\bar{b}$ and
$d(x,U*\nu)$ by $g(\underline{U}\langle\nu,\mathbf{1}\rangle)$
where $g(\zeta):=\inf_{x\in{\cal X}}d(x,\zeta)$ in the indicator
functions of each terms of the construction~(\ref{bpe}) of the
process $X^K_t$, we can stochastically dominate the population
size $\langle X^K_t,\mathbf{1}\rangle$ by a birth and death Markov
process $(Z^K_t)_{t\geq 0}$ with initial state $Z_0^K=\langle
  X^K_0,1\rangle$ and transition
rates
\begin{equation*}
  \begin{array}{ll}
    i2\bar{b}\quad & \mbox{from\ }i/K\mbox{\ to\ }(i+1)/K, \\
    ig(\underline{U}\frac{i}{K})\quad & \mbox{from\ }i/K\mbox{\ to\ }(i-1)/K.
  \end{array}
\end{equation*}
Therefore, it suffices to prove that $\sup_{K\geq 0}\sup_{t\geq
  0}E((Z^K_t)^q)<+\infty$.

Let us define $p^k_t={P}(Z^K_t=k/K)$. Then
\begin{align*}
  \frac{d}{dt}{E}((Z^K_t)^q) & =\sum_{k\geq
    1}\left(\frac{k}{K}\right)^q\frac{dp_t^k}{dt} \\
  & =\frac{1}{K^q}\sum_{k\geq 1}k^q\left[2\bar{b}(k-1)p^{k-1}_t
    +(k+1)g\left(\frac{k+1}{K}\right)p^{k+1}_t
    -k\left(2\bar{b}+g\left(\frac{k}{K}\right)\right)p^k_t\right] \\
  & =\frac{1}{K^q}\sum_{k\geq
    1}\left[2\bar{b}\left(\left(1+\frac{1}{k}\right)^q-1\right)
    +g\left(\frac{k}{K}\right)
    \left(\left(1-\frac{1}{k}\right)^q-1\right)\right]k^{q+1}p^k_t.
\end{align*}
Now, by~(H5)~(iii)~(\ref{eq:d-infty}),
$g(\alpha)\rightarrow+\infty$ when $\alpha\rightarrow +\infty$, so
there exists $\alpha_0$ such that, for any $\alpha\geq\alpha_0$,
$g(\alpha)\geq 4\bar{b}$. Therefore, for $k\geq K\alpha_0$,
$2\bar{b}((1+1/k)^q-1)+g(k)((1-1/k)^q-1)\leq
-2\bar{b}[3-2(1-1/k)^q-(1+1/k)^q]$, the RHS term being equivalent
to $-2\bar{b}q/k$. Therefore, enlarging $\alpha_0$ if necessary
and using the fact that $(1+\alpha)^q-1\leq \alpha(2^q-1)$ for any
$\alpha\in[0,1]$, we can write
\begin{align*}
  \frac{d}{dt}{E}((Z^K_t)^q) & \leq\sum_{k=1}^{\lceil
    K\alpha_0\rceil-1}2\bar{b}(2^q-1)\left(\frac{k}{K}\right)^q
  -\sum_{k\geq\lceil K\alpha_0\rceil}\bar{b}q\left(\frac{k}{K}\right)^qp^k_t \\
  & \leq 2\bar{b}(2^q-1)\alpha_0^2+\bar{b}q\alpha_0^2
  -\bar{b}q {E}((Z^K_t)^q).
\end{align*}
Writing $C=(2(2^q-1)+q)\alpha_0^2/q$, this differential inequality
solves as
\begin{equation*}
  {E}((Z^K_t)^q)\leq C+({E}(\langle X_0^K,1\rangle^q)-C)e^{-\bar{b}qt},
\end{equation*}
which gives the required uniform bound.\hfill$\Box$\\

\paragraph{Proof of Corollary~\ref{cor:IPS-TSS}}
Let $\Gamma$ be a measurable subset of ${\cal X}$. Let us prove
that
\begin{equation}
  \label{eq:pf-cor}
  \lim_{K\rightarrow+\infty}
   E\big[\langle X^K_{t/Ku_K},\mathbf{1}_{\Gamma}\rangle\big]
  = E\big[\bar{n}(Y_t)\mathbf{1}_{Y_t\in\Gamma}\big].
\end{equation}
By (H5)-(iii)-(\ref{eq:d-infty}), there exists $\zeta_0>0$ such that
for any $\zeta>\zeta_0$ and $x\in{\cal X}$, $d(x,\zeta)>\bar{b}$.
Therefore, by~(H5)-(iv), for any $x\in{\cal X}$,
$\bar{n}(x)\in[0,\zeta_0/\underline{U}]$. Fix $\varepsilon>0$, and
write $[0,\zeta_0/\underline{U}]\subset\cup_{i=1}^{p}I_i$, where $p$
is the integer part of $\zeta_0/(\underline{U}\varepsilon)$, and
$I_i=[(i-1)\varepsilon,i\varepsilon[$. Define $\Gamma_i=\{x\in{\cal
  X}:\bar{n}(x)\in I_i\}$ for $1\leq i\leq p$, and
apply~(\ref{eq:def-cvgce}) to the sets
$\Gamma\cap\Gamma_1,\ldots,\Gamma\cap\Gamma_p$ with $n=1$, $t_1=t$ and
the constant $\varepsilon$ above. Then, by Lemma~\ref{lem:moment}, for
some constant $C>0$ and for sufficiently large $K$,
\begin{align*}
  \limsup_{K\rightarrow+\infty}  E\big[\langle
  X^K_{t/Ku_K},\mathbf{1}_{\Gamma}\rangle\big] &
  \leq\limsup_{K\rightarrow+\infty} E\big[\langle
  X^K_{t/Ku_K},\mathbf{1}_{\Gamma}\rangle
  \mathbf{1}_{\langle X^K_{t/Ku_K},\mathbf{1}\rangle\leq
  C}\big]+\varepsilon \\
  & \leq\sum_{i=1}^p\limsup_{K\rightarrow+\infty} E\big[\langle
  X^K_{t/Ku_K},\mathbf{1}_{\Gamma\cap\Gamma_i}\rangle
  \mathbf{1}_{\langle X^K_{t/Ku_K},\mathbf{1}\rangle\leq C}\big]
  +\varepsilon \\
  & \leq\sum_{i=1}^p(i+1)\varepsilon P(Y_t\in\Gamma\cap\Gamma_i)
  +\varepsilon \\
  & \leq\sum_{i=1}^p\big(
   E\big[\bar{n}(Y_t)\mathbf{1}_{X_t\in\Gamma\cap\Gamma_i}\big]
  +2\varepsilon P(Y_t\in\Gamma_i)\big)+\varepsilon \\
  & \leq E\big[\bar{n}(Y_t)\mathbf{1}_{Y_t\in\Gamma}\big]+3\varepsilon.
\end{align*}
A similar estimate for the \emph{lim\:inf} ends the proof
of~(\ref{eq:pf-cor}), which implies the convergence of one-dimensional
laws for the required topology.

The same method gives easily the required limit when we consider a
finite number of times $t_1,\ldots,t_n$.\hfill$\Box$\\

Observe that the fact that the limit process is not
right-continuous prevents the possibility to obtain a convergence
for the Skorohod topology on $\mathbb{D}([0,T],M_F({\cal X}))$.

\subsection{Idea of the proof}
\label{sec:idea}

Theorem~\ref{thm:IPS-TSS} can be proved in a similar way as in
Champagnat~\cite{Ch04}. Let us give an idea of the method in order
to explain the assumptions, the various parameters appearing in
Theorem~\ref{thm:IPS-TSS} and the tools involved in the proof. It
is based on two ingredients: the study of a monomorphic population
before the first mutation, and the study of the invasion of a
single mutant individual in this population.

1) \:The first part obtains from large deviation results for the
convergence of $ X^K_t$ to $n_t(x)\delta_x$ when the initial
population is monomorphic with trait $x$, where $n_t(x)$
satisfies~(\ref{eq:monomorph}). Any positive solution
to~(\ref{eq:monomorph}) converges to $\bar{n}(x)$ when
$t\rightarrow+\infty$, and hence reaches a given neighborhood of
$\bar{n}(x)$ in finite time, i.e.\ on an infinitesimal time scale
with respect to the mutation time scale. Large deviations theory
allows us to show that the exit time of $\langle
X^K_t,\mathbf{1}\rangle\:$ from this neighborhood behaves as
$\exp(KC)$ for some $C>0$ (problem of exit from a domain, Freidlin
and Wentzell~\cite{FW84}). Thanks to the right part of
Assumption~(\ref{eq:K-u_K-bis}), we can prove that, with high
probability, $\langle X^K_t,\mathbf{1}\rangle$ is close to
$\bar{n}(x)$ when the first mutation occurs. Therefore, the total
mutation rate is close to
$u_K\mu(x)K\bar{n}(x)b(x,V(0)\bar{n}(x))$ and so, on the mutation
time scale $t/Ku_K$, the rate of mutation is close to
$\bar{n}(x)\mu(x)b(x,V(0)\bar{n}(x))$, which explain the left part
of the RHS of~(\ref{eq:g}). This argument can be made rigorous
using stochastic domination results similar to the one used at the
beginning of the proof of Lemma~\ref{lem:moment}, and leads to the
following result:

\begin{lem}
  \label{lem:tau}
  Let $\tau_1$ denote the first mutation time and
  $\mathbf{P}^K_{ X^K_0}$ the law of $ X^K$ with initial state
  $ X^K_0$.
    Given $x\in{\cal X}$ and a sequence of integers $(z_K)_{K\geq 1}$
    such that $z_K/K\rightarrow z>0$,
  \begin{description}
  \item[\textmd{(a)}] For any $\varepsilon>0$,
    \begin{equation}
      \label{eq:lem-tau-(b)}
      \lim_{K\rightarrow+\infty}\mathbf{P}^K_{\frac{z_K}{K}\delta_x}
      \left(\tau_1>\log K,\ \sup_{t\in[\log K,\tau_1]}|\langle
      X^K_t,\mathbf{1}\rangle-\bar{n}(x)|>\varepsilon\right)=0
    \end{equation}
    and
    \begin{equation*}
      \lim_{K\rightarrow +\infty}\mathbf{P}^K_{\frac{z_K}{K}\delta_x}
      (\tau_1<\log K)=0.
    \end{equation*}
    In particular, under $\mathbf{P}^K_{\frac{z_K}{K}\delta_x}$,
    $X^K_{\log K}\rightarrow\bar{n}(x)\delta_x$ and
    $X^K_{\tau_1-}\rightarrow\bar{n}(x)\delta_x$ in probability.

  \item[\textmd{(b)}] For any $t>0$,
    \begin{equation*}
      \lim_{K\rightarrow+\infty}\mathbf{P}^K_{\frac{z_K}{K}\delta_x}
      \left(\tau_1>\frac{t}{Ku_K}\right)
      =\exp\big(-\beta(x)t\big),
    \end{equation*}
    where $\beta(x)=\mu(x)\bar{n}(x)b(x,V(0)\bar{n}(x))$.
  \end{description}
\end{lem}

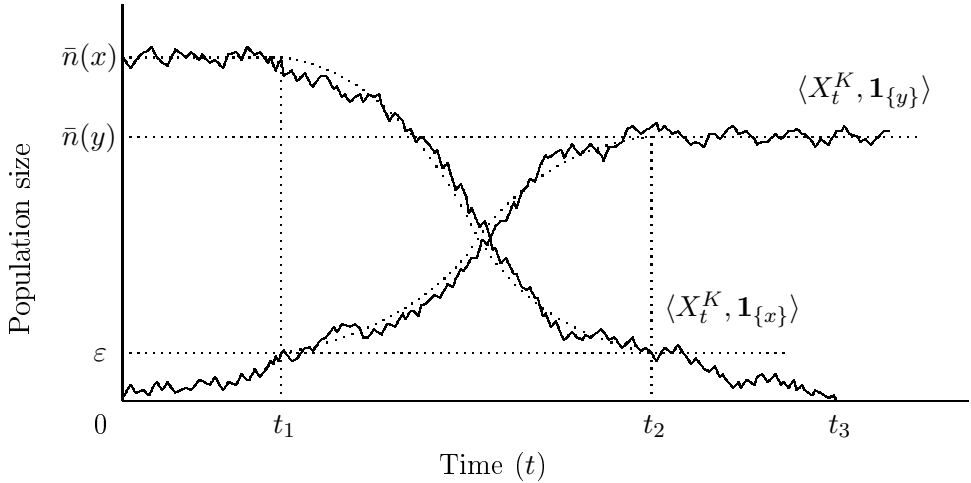
\begin{figure}[ht]
  \begin{center}
    \begin{picture}(380,190)(-45,-23)
      \put(0,0){\line(1,0){320}} \put(0,0){\line(0,1){150}}
      \put(-11,-12){0} \put(-11,15){$\varepsilon$}
      \dottedline{3}(0,18)(250,18)
      \put(-23,97){$\bar{n}(y)$} \put(-23,127){$\bar{n}(x)$}
      \dottedline{3}(0,100)(300,100) \dottedline{3}(0,130)(60,130)
      \put(-43,25){\rotatebox{90}{Population size}} \put(57,-12){$t_1$}
      \put(197,-12){$t_2$} \put(267,-12){$t_3$}
      \put(140,-25){\makebox(0,0){Time ($t$)}}
      \dottedline{3}(60,0)(60,130) \dottedline{3}(200,0)(200,100)
      \qbezier[27](60,130)(100,125)(130,70)
      \qbezier[25](130,70)(150,28)(200,18)
      \qbezier[25](60,18)(105,26)(130,55)
      \qbezier[25](130,55)(155,95)(200,100)
      \dottedline{0.5}(0,1)(3,6)(5,3)(8,4)(10,2)(12,6)(13,7)(15,4)(17,7)%
(20,3)(24,5)(26,4)(29,8)(33,6)(35,9)(38,4)(40,7)(43,6)(45,10)(47,8)(48,8)%
(50,12)(52,9)(55,13)(56,15)(57,14)(58,17)(59,15)(60,18)(61,19)(64,15)%
(66,16)(69,20)(71,20)(72,19)(75,23)(77,24)(80,28)(81,26)(83,29)(86,28)%
(88,29)(90,27)(91,24)(93,24)(95,27)(98,25)(100,29)(101,29)(103,28)(106,33)%
(107,31)(109,34)(110,34)(112,32)(113,36)(115,37)(118,40)(120,39)(123,44)%
(125,43)(126,46)(128,50)(129,49)(132,54)(134,54)(136,61)(138,59)(141,66)%
(143,68)(145,68)(148,75)(149,73)(151,79)(153,81)(154,84)(156,83)(158,89)%
(160,93)(161,92)(164,94)(166,95)(167,92)(170,97)(173,93)(175,95)(177,97)%
(181,93)(182,91)(184,96)(186,96)(189,98)(191,103)(193,104)(196,101)(199,103)%
(202,105)(203,102)(205,104)(208,101)(210,102)(213,99)(215,97)(217,98)%
(220,96)(222,101)(224,101)(226,99)(229,103)(231,104)(233,101)(236,102)%
(239,99)(240,97)(242,98)(244,98)(247,101)(249,99)(252,102)(254,102)(255,103)%
(257,101)(258,101)(261,98)(262,96)(264,99)(267,97)(270,102)(272,102)%
(273,104)(275,101)(278,100)(279,98)(281,99)(283,97)(288,102)(290,102)
      \put(255,115){$\langle X^K_t,\mathbf{1}_{\{y\}}\rangle$}
      \dottedline{0.5}(0,130)(3,126)(6,130)(8,132)(9,132)(11,134)(13,130)%
(15,128)(17,127)(20,133)(24,130)(26,132)(29,129)(30,131)(33,128)(35,129)%
(38,126)(40,131)(43,132)(45,129)(47,134)(48,133)(50,131)(52,132)(55,128)%
(56,130)(57,127)(58,125)(59,128)(60,129)(61,123)(64,125)(66,121)(69,124)%
(71,120)(72,119)(75,123)(77,123)(80,118)(81,119)(83,115)(84,116)(86,113)%
(88,115)(90,116)(91,113)(93,116)(95,115)(97,116)(98,112)(100,108)(101,106)%
(103,109)(106,106)(107,102)(109,103)(111,98)(113,100)(116,97)(118,94)%
(119,95)(121,88)(123,90)(124,86)(126,80)(128,82)(130,77)(131,71)(132,73)%
(135,66)(137,68)(139,62)(140,57)(141,58)(143,52)(145,54)(148,48)(149,48)%
(151,43)(153,42)(154,43)(156,39)(157,39)(159,32)(160,33)(162,28)(164,30)%
(166,25)(167,23)(170,26)(172,24)(174,23)(175,25)(177,24)(181,27)(182,26)%
(184,23)(186,24)(189,20)(191,22)(193,24)(194,21)(196,19)(198,20)(200,18)%
(202,15)(203,17)(205,20)(208,18)(210,21)(213,17)(215,14)(217,16)(220,12)%
(222,11)(224,13)(226,9)(227,10)(229,7)(231,5)(233,7)(234,5)(236,8)%
(239,6)(240,9)(242,10)(244,7)(247,10)(249,6)(250,7)(252,5)(254,8)%
(255,4)(257,6)(258,4)(261,2)(262,4)(264,3)(267,1)(269,2)(270,0)
      \put(205,32){$\langle X^K_t,\mathbf{1}_{\{x\}}\rangle$}
    \end{picture}
  \end{center}
  \caption{The three steps of the invasion and fixation of a
    mutant trait $y$ in a monomorphic population with trait $x$.
    Plain curves represent the resident and mutant densities $\langle
    X^K_t,\mathbf{1}_{\{x\}}\rangle$ and $\langle
    X^K_t,\mathbf{1}_{\{y\}}\rangle$, respectively. Dotted curves
    represent the solution of Eq.~(\ref{eq:dimorph}) with initial
    state $n_0(x)=\bar{n}(x)$ and $n_0(y)=\varepsilon$.}
  \label{fig:inv-fix}
\end{figure}

2) \:The study of the invasion of a mutant individual with trait
$y$ can be divided in three steps represented in
Fig.~\ref{fig:inv-fix}.

Firstly, the invasion of the mutant (between 0 and $t_1$ in
Fig.~\ref{fig:inv-fix}) can be defined as the growth of the mutant
density $\langle X^K_t,\mathbf{1}_{\{y\}}\rangle$ from $1/K$ (one
individual) to a fixed small level $\varepsilon$ ($\varepsilon K$
individuals). As long as the mutant density is small, the dynamics
of the resident density $\langle X^K_t,\mathbf{1}_{\{x\}}\rangle$
is close to the one it followed before the mutation, so it is
close to $\bar{n}(x)$ with high probability. Therefore, between
$0$ and $t_1$, the birth and death rates of an individual with
trait $y$ are close to $b(y,V(y-x)\bar{n}(x))$ and
$d(y,U(y-x)\bar{n}(x))$ respectively. Therefore, the number of
mutant individuals is close to a binary branching process with the
parameters above. When $K\rightarrow +\infty$, the probability
that such a branching process reaches level $\varepsilon K$ is
close to its survival probability, which writes
$[f(y,x)]_+/b(y,V(y-x)\bar{n}(x))$. This gives the second part of
the RHS of~(\ref{eq:g}).

Secondly, once the invasion succeeded (which is possible only if
$f(y,x)>0$), the dynamics of the densities of traits $x$ and $y$
are close to the solution to the dimorphic logistic
equation~(\ref{eq:dimorph}) with initial state
$(\bar{n}(x),\varepsilon)$, represented in dotted curves between
$t_1$ and $t_2$ in Fig.~\ref{fig:inv-fix}. Because of
Assumption~(H6), the resident density can be proved to reach level
$\varepsilon$ with high probability (at time $t_2$ in
Fig.~\ref{fig:inv-fix}).

Finally, a similar argument as in the first step above allows us
to prove that the resident population density $\langle
X^K_t,\mathbf{1}_{\{x\}}\rangle$ follows approximately a binary
branching process with birth rate $b(y,V(x-y)\bar{n}(y))$ and
death rate $d(y,U(x-y)\bar{n}(y))$. Since $f(x,y)<0$ by
Assumption~(H6), this is a sub-critical branching process, and
therefore, the resident trait $x$ disappears in finite time $t_3$
with high probability.

We can show, using results on branching processes, that $t_1$ and
$t_3-t_2$ are of order $\log K$, whereas $t_2-t_1$ depends only on
$\varepsilon$. Therefore, the left part of~(\ref{eq:K-u_K-bis})
ensures that the three steps of the invasion are completed before
the next mutation, with high probability. The previous heuristics
can be made rigorous using further comparison results, and leads
to the following result.

\begin{lem}
  \label{lem:theta}
  Assume that the initial population is made of individuals with
  traits $x$ and $y$ satisfying assumption~(H6)~(i) or~(ii).
  Let $\theta_0$ denote the first time when the population gets monomorphic,
  and $V_0$ the remaining trait. Let
  $(z_K)_{K\geq 1}$ be a sequence of integers such that
  $z_K/K\rightarrow\bar{n}(x)$. Then,
  \begin{gather}
    \label{eq:lem-theta-1}
    \lim_{K\rightarrow+\infty}\mathbf{P}^K_{\frac{z_K}{K}\delta_x
      +\frac{1}{K}\delta_y}(V_0=y)
    =\frac{[f(y,x)]_+}{b(y,V(y-x)\bar{n}(x))}, \\
    \label{eq:lem-theta-2}
    \lim_{K\rightarrow+\infty}\mathbf{P}^K_{\frac{z_K}{K}\delta_x
      +\frac{1}{K}\delta_y}(V_0=x)
    =1-\frac{[f(y,x)]_+}{b(y,V(y-x)\bar{n}(x))}, \\
    \label{eq:lem-theta-3}
    \forall\eta>0,\quad
    \lim_{K\rightarrow+\infty}\mathbf{P}^K_{\frac{z_K}{K}\delta_x
      +\frac{1}{K}\delta_y}
    \left(\theta_0>\frac{\eta}{Ku_K}\wedge\tau_1\right)=0 \\
    \label{eq:lem-theta-4}
    \mbox{and}\quad\forall\varepsilon>0,\quad
    \lim_{K\rightarrow+\infty}\mathbf{P}^K_{\frac{z_K}{K}\delta_x
      +\frac{1}{K}\delta_y}\left(|\langle X^K_{\theta_0},
      \mathbf{1}\rangle-\bar{n}(V_0)|<\varepsilon\right)=1,
  \end{gather}
  where $f(y,x)$ has been defined in~(\ref{eq:def-fitness}).
\end{lem}

Once these lemmas are proved, the proof can be completed by
observing that the generator $A$ of the process $(Y_t,t\geq 0)$ of
Theorem~\ref{thm:IPS-TSS} can be written as
\begin{equation}
  \label{eq:generator-TSS-mu(x,dh)}
  A\varphi(x)=\int_{\mathbb{R}^l}(\varphi(y)-\varphi(x))
  \beta(x)\kappa(x,dy),
\end{equation}
where $\beta(x)$ has been defined in Lemma~\ref{lem:tau} and the
probability measure $\kappa(x,dh)$ is defined by
\begin{multline}
  \label{eq:def-kappa}
  \kappa(x,dy)=\left(1-\int_{\mathbb{R}^l}
    \frac{[f(z,x)]_+}{b(z,V(z-x)\bar{n}(x))}M(x,z)dz\right)
  \delta_x(dy) \\ +\frac{[f(y,x)]_+}{b(y,V(y-x)\bar{n}(x))}M(x,y)dy.
\end{multline}
This means that the process $Y$ with initial state $x$ can be
constructed as follows: let $(M(k),k=0,1,2,\ldots)$ be a Markov
chain in ${\cal X}$ with initial state $x$ and with transition
kernel $\kappa(x,dy)$, and let $(N(t),t\geq 0)$ be an independent
standard Poisson process. Let also $(T_n)_{n\geq 1}$ denote the
sequence of jump times of the Poisson process $N$. Then, the
process $(Y_t,t\geq 0)$ defined by
\begin{equation*}
  Y_t:=M\left(N\left(\int_0^t\beta(Y_s)ds\right)\right)
\end{equation*}
is a Markov process with infinitesimal
generator~(\ref{eq:generator-TSS-mu(x,dh)})
(cf.~\cite{Ethier-Kurtz} chapter~6).

Let $P_x$ denote its law, and define $(S_n)_{n\geq 1}$ by
$T_n=\int_0^{S_n}\beta(Y_s)ds$. Observe that any jump of the
process $Y$ occurs at some time $S_n$, but that all $S_n$ may not
be effective jump times for $Y$, because of the Dirac mass at $x$
appearing in~(\ref{eq:def-kappa}).

Fix $t>0$, $x\in{\cal X}$ and a measurable subset $\Gamma$ of
${\cal
  X}$. Under $P_x$, $S_1$ and $Y_{S_1}$ are independent,
$S_1$ is an exponential random variable with parameter $\beta(x)$,
and $Y_{S_1}$ has law $\kappa(x,\cdot)$.  Therefore, for any
$n\geq 1$, the strong Markov property applied to $Y$ at time $S_1$
yields
\begin{multline}
  \label{eq:induction-X}
  P_x(S_n\leq t<S_{n+1},\ Y_t\in\Gamma) \\
  =\int_0^t\beta(x)e^{-\beta(x)s}\int_{\mathbb{R}^l}
  \mathbf{P}_{y}(S_{n-1}\leq t-s<S_n,\ Y_{t-s}\in\Gamma)
  \kappa(x,dy)ds
\end{multline}
and
\begin{equation}
  \label{eq:init-induction-X}
  P_x(0\leq t<S_1,\ Y_t\in\Gamma)
  =\mathbf{1}_{\{x\in\Gamma\}}e^{-\beta(x)t}.
\end{equation}

Using the Markov property at time $\tau_1$ and
Lemmas~\ref{lem:tau} and~\ref{lem:theta}, we can prove that, when
we replace $S_n$ by the $n$-th mutation time of $X^K_{t/Ku_K}$ and
$Y_t$ by the support of $X^K_{t/Ku_K}$ (when it is a singleton) in
the LHS of~(\ref{eq:induction-X}) and~(\ref{eq:init-induction-X}),
the same relations hold in the limit $K\rightarrow +\infty$.
Therefore, Theorem~\ref{thm:IPS-TSS} is proved for one-dimensional
time marginals. A similar method generalizes to finite dimensional
laws.

\end{document}